\documentclass[11pt]{article}

\usepackage{amssymb,latexsym, amsmath}
\usepackage{graphicx}

\newtheorem{theorem}{Theorem}[section]
\newtheorem{coro}[theorem]{Corollary}
\newtheorem{defi}[theorem]{Definition}
\newtheorem{lemma}[theorem]{Lemma}
\newtheorem{prop}[theorem]{Proposition}

\def\slfrac#1#2{\hbox{\kern.1em %
 \raise.5ex\hbox{\the\scriptfont0 #1}\kern-.11em %
 /\kern-.15em\lower.25ex\hbox{\the\scriptfont0 #2}}}

\newcommand{\eqn}[1]{(\ref{#1})}

\newcommand{\eeq}{\end{equation}}
\newcommand{\beql}[1]{\begin{equation}\label{#1}}

\newcommand{\loc}{loc}

\newcommand{\ZZ}{{\mathbb Z}}
\newcommand{\RR}{{\mathbb R}}

\newcommand{\NN}{{\mathbb N}}

\newcommand{\QQ}{{\mathbb Q}}

\newcommand{\dM}{{N}^{1}}
\newcommand{\dN}{{D}}

\newcommand{\sB}{{\cal B}}
\newcommand{\sBp}{{\cal B}^{'}}

\newcommand{\sBo}{{B_{\emptyset}}}
\newcommand{\sC}{{\cal C}}

\newcommand{\sG}{{\cal G}}

\newcommand{\sL}{{\cal L}}

\newcommand{\tauL}{\tau^{L}}
\newcommand{\tauS}{\tau^{S}}
\newcommand{\tb}{\bar{b}}
\DeclareMathOperator{\var}{Var}
\newcommand{\ddL}{\Omega^{L}}

\newcommand{\muT}{{\mu_S}}
\DeclareMathOperator{\meas}{meas}
\DeclareMathOperator{\Prob}{Prob}

\renewcommand{\ll}{\langle\!\langle}
\renewcommand{\gg}{\rangle\!\rangle}

\makeatletter
\def\@sect#1#2#3#4#5#6[#7]#8{\ifnum #2>\c@secnumdepth
     \def\@svsec{}\else
     \refstepcounter{#1}\edef\@svsec{\csname the#1\endcsname.\hskip .75em }\fi
     \@tempskipa #5\relax
      \ifdim \@tempskipa>\z@
        \begingroup #6\relax
          \@hangfrom{\hskip #3\relax\@svsec}{\interlinepenalty \@M #8\par}%
        \endgroup
       \csname #1mark\endcsname{#7}\addcontentsline
         {toc}{#1}{\ifnum #2>\c@secnumdepth \else
                      \protect\numberline{\csname the#1\endcsname}\fi
                    #7}\else
        \def\@svsechd{#6\hskip #3\@svsec #8\csname #1mark\endcsname
                      {#7}\addcontentsline
                           {toc}{#1}{\ifnum #2>\c@secnumdepth \else
                             \protect\numberline{\csname the#1\endcsname}\fi
                       #7}}\fi
     \@xsect{#5}}
\def\@begintheorem#1#2{\it \trivlist \item[\hskip \labelsep{\bf #1\ #2.}]}

\def\plain{plain}\ifx\fmtname\plain\csname fi\endcsname
     
     \input docstrip
     \preamble

     Do not distribute the stripped version of this file.
     The checksum in the header refers to the documented version.

     \endpreamble
     \generateFile{here.sty}{t}{\from{here.doc}{}}
     \endinput
\fi
\ifcat a\noexpand @\let\next\relax\else\def\next{%
    \documentstyle[here,doc]{article}\MakePercentIgnore}\fi\next
\ifx\@Hxfloat\@Hundef\else\expandafter\endinput\fi
\let\@Hxfloat\@xfloat
\def\@xfloat#1[{\@ifnextchar{H}{\@HHfloat{#1}[}{\@Hxfloat{#1}[}}
\def\@HHfloat#1[H]{%
\expandafter\let\csname end#1\endcsname\end@Hfloat
\vskip\intextsep\vbox\bgroup\def\@captype{#1}\parindent\z@
\ignorespaces}
\def\end@Hfloat{\egroup\vskip \intextsep}

\makeatother

\begin{document}

\begin{center}
{\Large 
{\bf   Level Sets of the Takagi Function: Local Level Sets} 
}\\

\vspace{1.5\baselineskip}
{\em Jeffrey C. Lagarias
\footnote{This author's work was supported by NSF Grant DMS-0801029 and DMS-1101373.} }\\
\vspace*{.2\baselineskip}
Dept. of Mathematics \\
University of Michigan \\
Ann Arbor, MI 48109-1043\\
\vspace*{1.5\baselineskip}
{\em Zachary Maddock}
\footnote{This author's work was supported by the NSF
  through a Graduate Research Fellowship.} \\
\vspace*{.2\baselineskip}
Dept. of Mathematics \\
Columbia University  \\
New York, NY  10027\\

\vspace*{2\baselineskip}
(March 2, 2012) \\
\vspace{3\baselineskip}
{\bf ABSTRACT}
\end{center}
The Takagi function $\tau: [0,1] \to [0, 1]$ is a 
continuous non-differentiable function
constructed by Takagi in 1903.
The
 level sets $L(y)=\{ x: \tau(x)=y \}$ 
 of the Takagi function $\tau(x)$ are studied by introducing a notion of 
 local level set into which level sets are partitioned. 
 Local level sets are simple to analyze,
 reducing questions to understanding the relation of level
 sets to local level sets, which is more complicated.
  It is known  that for a ``generic"  full  Lebesgue measure set of ordinates $y$,
 the level sets are  finite sets.
 Here it is shown for a ``generic" 
 full Lebesgue measure set of abscissas $x$, the
 level set $L(\tau(x))$  is uncountable.
 An interesting singular monotone function is constructed associated
 to local level sets, and is used to show the expected number
 of local level sets at a random level $y$ is exactly $\frac{3}{2}$.\\

\noindent



%
%
%
%
\setlength{\baselineskip}{1.0\baselineskip}

\section{Introduction}

The Takagi function $\tau(x)$ is  a function
defined on the unit interval $ x\in [0,1]$
which was introduced by Takagi \cite{Takagi} in 1903
as an example of a continuous nondifferentiable function.
It can be defined by 
\begin{equation}\label{eq101}
\tau(x) := \sum_{n=0}^{\infty} \frac{\ll 2^n x \gg}{2^n}
\end{equation}
where $\ll x \gg := \inf_{n \in \ZZ} |x-n|$ is the distance from $x$ to the nearest integer.
Variants of this function were presented by van der Waerden \cite{vdW30} in
1930 and de Rham \cite{deR57} in 1957.

%
%

\begin{figure}[h]
\centering
\includegraphics[width=3in]{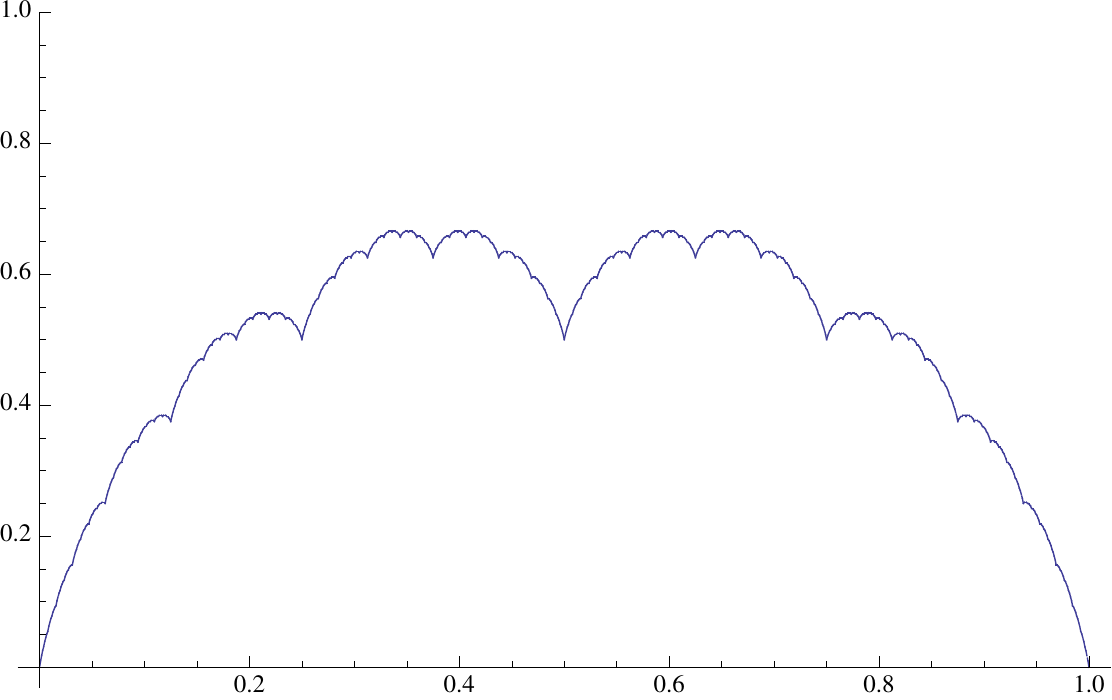}
\caption{Graph of the Takagi function $\tau(x)$.}
\label{fig11}
\end{figure}

An alternate interpretation of the Takagi function involves the symmetric tent map
$T: [0,1] \to [0,1]$, given by 
\begin{equation}
\label{eq101a}
T(x) = \left\{ 
\begin{array}{cc} 
2x & \mbox{if}~~ 0\le x \le \frac{1}{2},\\
~&~\\
2-2x & \mbox{if}~~\frac{1}{2}\le x \le 1
\end{array}
\right.
\end{equation}
(see  \cite{FL00} for  further references). 
Then we have
\begin{equation*}
\tau(x) := \frac{1}{2} \left( \sum_{n=1}^{\infty} \frac{1}{2^n} T^{(n)}(x) \right),
\end{equation*}
where $T^{(n)}(x)$ denotes the $n$-th iterate of $T(x)$.
The Takagi function appears in many contexts and has been studied
extensively;  see  the recent surveys of   Allaart and Kawamura \cite{AK11}
 and of the first author \cite{La11}.

In this paper we consider certain properties of the graph of the Takagi function
$$
\sG(\tau):= \{ (x, \tau(x)): ~0 \le x \le 1 \},
$$
which is pictured in Figure~\ref{fig11}. 
It is well known that the values of the Takagi function satisfy $0 \le \tau(x) \le \frac{2}{3}.$
It is also known that this graph has Hausdorff dimension $1$ in $\RR^2$, see Mauldin
and Williams \cite[Theorem 7]{MW86}), and furthermore it is
$\sigma$-finite, see Anderson and Pitt \cite[Thm. 6.4]{AP89}. 
Here we  study  the structure of the level sets of this graph. We make the
following definition, which contains a special convention concerning
dyadic rationals which simplifies theorem statements. 

%
%

\begin{defi}~\label{de11}
{\em 
For $0 \le y \le \frac{2}{3}$  the 
{\em (global) level set}  $L(y)$ at level $y$ is 
$$
L(y):= \{ x: ~ \tau(x) =y, ~~0 \le x \le 1 \}. 
$$
We make the convention that 
$x$  specifies a binary expansion; thus each dyadic rational
value  $x= \frac{m}{2^n}$  in a level set occurs twice, 
labeled by  its two possible
 binary expansions. 
(Technically $L(y)$ is a multiset, with multiplicities $1$ or $2$.)
}
\end{defi}

Level sets have a complicated structure, depending on the value of $y$.
It is known that there are different levels $y$ where the level set  $L(y)$ 
is finite, countably infinite, or
uncountably infinite, respectively.
In 1959 Kahane \cite[Sec. 1]{Ka59} noted that 
the level set $L(\frac{2}{3})$ 
 was a perfect, totally disconnected set of
Lebesgue measure $0$, and in  1984
Baba \cite{Baba84} showed that $L(\frac{2}{3})$ has Hausdorff dimension
$\frac{1}{2}$. 
The second author recently proved (\cite{Maddock09}) that the Hausdorff
dimension of any level set is at most $0.668$, and conjectured that
the example of Baba achieves the largest possible Hausdorff dimension.  
Recently de Amo, Bhouri, D\'{i}az Carrillo and Fern\'{a}ndez-S\'{a}nchez \cite{ABDF11} 
proved this conjecture.
The level sets at a rational level $y \in \QQ$ are particularly interesting.
 Knuth \cite[Sect. 7.2.1.3, Exercises 82-85]{Kn09} gave a (not necessarily halting)
 algorithm 
 to determine the structure of level
sets at rational levels, revealing very complicated behaviors.  For example he determined that
$L(\frac{1}{5})$ has a finite level set with exactly two elements, namely
$~L(\frac{1}{5}) = \{ x, 1-x\}~~\mbox{with}~~ x= \frac{83581}{87040}.$ 
He also noted that
$L(\frac{1}{2})$ is countably infinite; we include a proof in Theorem~\ref{th38} below.
In 2008 Buczolich \cite{Buz08} 
proved that, in the sense of Lebesgue measure on $y \in [0, \frac{2}{3}]$, 
almost all level sets $L(y)$ are finite sets. \\

The object of this paper is to introduce and study 
 the notion of ``local level set". These are sets determined
locally by combinatorial operations on the binary expansion of a real number $x$; 
they are closed sets and we 
show that  each level set  decomposes into a disjoint union of local level sets. 
(The convention on dyadic rationals made in the definition 
above is needed for disjointness of the union.)
The structure
of local level sets is completely analyzable: they are either finite sets or Cantor
sets. Information about the  Hausdorff dimension of such sets can
readily be deduced from properties of the binary expansion of $x$. \\

We then study the relation of local level sets and level sets. 
How many local level sets are there in a given level set? 
To approach this question, we  study  the behavior of the Takagi function
restricted to the set $\ddL$ of left hand (abscissa) endpoints $x$
of all the local level sets; these endpoints parametrize the totality of all local level sets. 
We show that $\ddL$  is a closed perfect set (Cantor set) which has Lebesgue measure
$0$;  in a sequel (\cite{LM10b})
we show it has Hausdorff dimension $1$. We  show the Takagi function behaves relatively
nicely when restricted to $\ddL$, namely that $\tauS(x):=\tau(x)+ x$ is a 
monotone singular continuous function on this set.
It is therefore the integral of a singular probability measure on $[0,1]$,
which we call the Takagi singular measure.
Using this function we deduce that the expected
number of local level sets at a random level $0 \le y \le \frac{2}{3}$ is finite, and we
determine that  this expected value 
is exactly $\frac{3}{2}$. We also show that there is a dense set of values $y$ having
an infinite number of distinct local level sets.  

\medskip

Local level sets provide a way to take apart level sets and better understand their
structure. 
 In the rest of the introduction we state the main results of this paper in more detail.


\subsection{Local level sets}

This notion of local level set is attached to the binary expansion of abscissa point $x \in [0,1]$.
We show that certain combinatorial flipping operations applied
to the binary expansion of $x$ yield new points $x'$ in the same level set.
The totality of points reachable from $x$ by these combinatorial operations
will comprise the local level set $L_x^{loc}$ associated to $x$.

To describe this, let   $x \in [0,1]$ have a binary expansion:
\begin{equation*}
x := \sum_{j=1}^{\infty} \frac{b_j}{2^j}= 0.b_1 b_2 b_3..., ~~~~ \mbox{each}~b_j \in \{0, 1\}. 
\end{equation*}
The {\em flip operation} (or {\em complementing operation}) on  a single binary digit $b$ is 
\begin{equation*}
\bar{b} := 1- b.
\end{equation*}
We associate to the binary expansion the {\em digit sum function} $\dM(x)$ given by
\begin{equation*}
\dM_j(x) := b_1 +b_2 + \cdots + b_j.
\end{equation*}
We also associate to the binary expansion the  {\em  deficient digit function} $\dN_j(x)$ given by
\begin{equation*}
\dN_j(x):= j- 2\dM_j(x)  = j-  2(b_1+b_2+ \cdots + b_j) .
\end{equation*}
Here  $\dN_j(x)$ counts the excess of binary digits $b_k=0$ over those with $b_k=1$
in the first $j$ digits, i.e. it is positive if there are more $0$'s than $1$'s in the first $j$ digits.  
Note that for dyadic rationals 
$x= \frac{m}{2^n}$
the function values depend on which binary expansion  is used. \\

 We next associate  to any $x$  the sequence  of  digit positions $j$
at which tie-values $\dN_j(x)=0$ occur, which we call {\em balance points}; note that all such
$j$ are even.  The  {\em balance-set}  $Z(x)$ associated to $x$ is the set of balance
points, and is denoted
\beql{123}
Z(x) : = \{ c_k:~~\dN_{c_k}(x)=0\}.
\eeq
where we define $c_0=c_0(x) = 0$ and set $c_0(x)< c_1(x)< c_2(x) < ...$. This sequence of tie-values may
be finite or infinite. If it is finite, ending in $c_{n}(x)$, we make the
convention to adjoin a final ``balance point"  $c_{n+1}(x)= +\infty$.
 We call a {\em ``block"}  an indexed set of digits between two consecutive balance points, 
\begin{equation*}
B_k(x) := \{ b_j: ~c_k(x)  < j \le c_{k+1}(x)\},
\end{equation*}
which includes the second balance point but not the first.
We define an equivalence relation on
blocks, written  $B_k(x) \sim B_{k'}(x')$ to mean the block endpoints agree
($c_k(x)= c_{k'}(x')$ and $c_{k+1}(x) = c_{k'+1}(x')$)
and either $B_k(x) = B_{k'}(x')$ or $B_k(x) = \bar{B}_{k'}(x')$, where the bar operation flips
all the digits in the block, i.e. 
\begin{equation*}
b_j \mapsto \bar{b}_j:= 1- b_j,~~~~~~~ c_k < j \le c_{k+1}.
\end{equation*}
Finally, we define an equivalence relation $x \sim x'$ to mean that
they have identical balance-sets 
$Z(x) \equiv Z(x')$, and furthermore every block $B_k(x) \sim B_k(x')$ for $k \ge 0$.
Note that $x \sim 1-x$; this corresponds to a flipping operation being applied to
every binary digit.
We will show (Theorem \ref{th31})  
that the  equivalence relation $x\sim x'$ implies that $\tau(x)= \tau(x')$ so that $x$ and $x'$
are in the same level set of the Takagi function. 

%
%

\begin{defi}~\label{de12}
{\em 
The {\em local level set} $L_x^{loc}$ associated to $x$ is the set of equivalent points, 
\begin{equation*}
L_x^{\loc} := \{ x': ~~x' \sim x\}.
\end{equation*}
We use again the convention that $x$ and $x'$ denote binary
expansions, and hence dyadic rational numbers are represented by two
distinct binary expansions.}
\end{defi}

Each local level set  $L_x^{\loc}$  is 
a closed set. It  is  a finite set if the balance-set $Z(x)$ is finite, and  is a 
Cantor set (perfect
totally disconnected set)  if $Z(x)$ is infinite. 
Note that if $x$ is the expansion of 
a dyadic rational, then $L_x^{\loc}$ is finite and consists entirely of (expansions of) dyadic
rationals.

%
%
%
\begin{theorem}~\label{th14} {\em (Local level set partition)} \\
$~~~~~$ (1) Each local level set  $L_x^{\loc}$ is a  closed set contained
 in some level set.   
 
 (2) Two local level sets  $L_x^{loc}$
 and $L_{x'}^{loc}$ either coincide or are  disjoint.  
Thus each level set $L(y)$ partitions into a disjoint union of local level sets.
\end{theorem}

This easy result is  proved as part of Theorem \ref{th31} in Sect.~\ref{sec31}.
A priori this 
disjoint union could be finite, countable or uncountable. 
In Theorem \ref{th38}
we give an example of a level set that is a countably infinite  union
of local level sets; for this case $y= \frac{1}{2}$
is a dyadic rational. \smallskip

The Hausdorff dimension of a local level set $L_x^{loc}$ is 
restricted by the nature of its balance-set $Z(x)$. A necessary condition
to have positive Hausdorff dimension is that $Z(x)$  must have
positive upper asymptotic density in $\NN$.
This allows us to deduce the
following result.

%
%
%
\begin{theorem}~\label{th15} {\em (Generic local level sets)}
For a full Lebesgue measure set of abscissa points $x$  
the local level set $L_x^{loc}$ is a Cantor set (closed totally
disconnected perfect set) having  Hausdorff dimension $0$.
\end{theorem}

Theorem~\ref{th15} is proved in Sect.~\ref{sec32}.
This result implies that if an abscissa value $x$ is picked at random in $[0,1]$,
then with probability one the level set $L(\tau(x))$ is uncountably infinite.
This result differs strikingly from that of Buczolich \cite{Buz08},
who showed that if an ordinate value $y$ is picked uniformly in $[0, \frac{2}{3}]$
then the  level set $L(y)$   is finite with probability one.
 There is no inherent contradiction here:  drawing an abscissa value $x$
will preferentially select levels whose level set $L(\tau(x))$ is ``large", 
and Theorem~\ref{th15} quantifies  ``large."\\

In Sect ~\ref{sec33}  we  completely determine the structure of local level sets that contain
a rational number. We  prove they are either a finite set or a Cantor set of
positive Hausdorff dimension, and characterize when each case occurs (Theorem~\ref{th50}). 
One can  check directly that $x_0= \frac{1}{3}$ has a local level set of Hausdorff
dimension $ \frac{1}{2}$ (at level $y= \frac{2}{3}$), which shows that the Hausdorff dimension
upper bound of $\frac{1}{2}$ for level sets
obtained in \cite{ABDF11}  is sharp also for local level sets.


\subsection{Expected number of local level sets} 

Our second object is to relate local level sets to level sets. How many
local level sets belong to a given level set? 
To approach this problem   we first study (in Sect.~\ref{sec4}) the set $\ddL$ of all left
hand endpoints of local level sets, because this  set parameterizes all the local
level sets.  In Theorem \ref{th33} we 
 characterize its members in terms of their binary
expansions: they are exactly the values $x$ with binary expansions such that
$$
\dN_j(x) \ge 0 ~~~\mbox{for all} ~~~ j \ge 1.
$$
We call this latter  set the {\em deficient digit set}
and show it  is a closed, perfect set (Cantor set) of Lebesgue measure
zero.

In \S5 we define a new function, the {\em flattened Takagi function} $\tauL(x)$, 
which agrees with $\tau(x)$ on $\ddL$ and is defined
by linear interpolation across the gaps removed in constructing $\ddL$.
We prove  $\tauL(x)$ to be of bounded variation, and determine a Jordan decomposition.
  This consists of a nondecreasing piece $F_{+}(x) := \tauL(x) + x$
which we establish is a singular function whose points 
of increase are supported on $\ddL$,  and a 
strictly decreasing piece $F_{-}(x) := -x$ which is  absolutely continuous.
We name the function
$$
\tauS(x):= F_{+}(x) = \tauL(x) + x
$$
 the {\em Takagi singular function}, based on the following result,
 which is needed in the proof that the flattened Takagi function has
 bounded variation.

%
%
%
\begin{theorem} \label{th15a} {\em (Takagi singular function)}
The function  $\tauS(x)$ defined by $\tauS(x)= \tau(x) + x$ for
 $x \in \ddL$  
is a nondecreasing function on $\ddL$. Define its extension to
all $x \in [0,1]$ by 
$$
\tauS(x) := \sup\{ \tauS(x_1):  x_1 \le x ~~~\mbox{with}  ~~ x_1 \in \ddL\}.
$$
 Then the  function $\tauS(x)$  is  a monotone singular function. That is, it is
 a nondecreasing continuous function  
 having  
 $\tauS(0)=0, \tauS(1)=1$, which has derivative zero at (Lebesgue) almost all points of $[0,1]$. 
  The closure of the set of points of increase of $\tauS(x)$ is the deficient digit set $\ddL$.
\end{theorem}

In a sequel  (\cite{LM10b}) we study
 a nonnegative Radon measure $d\muT$, called the {\em Takagi singular measure},
such that 
\begin{equation*}
\tauS(x) = \int_{0}^x d\muT,
\end{equation*}
which is a probability measure on $[0,1]$.
This measure is singular with respect to Lebesgue measure.
There we show that
its support $\mbox{Supp}(\muT) = \ddL$ has (full) Hausdorff dimension $1$.
The Takagi singular  measure is not translation-invariant,
but it has  certain self-similarity properties under dyadic
rescalings. These are useful in explicitly computing
the measure of various interesting subsets of $\ddL$.
One may compare  analogous properties
of the Cantor function, see Dovghoshey et al \cite[Sect. 5]{DMRV06}.
\smallskip

The bounded variation property of the flattened Takagi function is used
to count the average number of local level sets, as follows. 

%
%
%
\begin{theorem}~\label{th16} {\em  (Expected number of local level sets)}
With respect to uniform (Lebesgue) measure on
the ordinate space $[0, \frac{2}{3}]$ a full measure set of points have a
finite number of local level sets. Furthermore
the  expected number of local level sets on a given level $y \in [0, \frac{2}{3}]$
is $\frac{3}{2}.$ 
\end{theorem}

This result is proved as Theorem~\ref{th37}, using the coarea formula for
functions of bounded variation.
We show that this result is non-trivial in that there are infinitely many levels
containing infinitely many distinct local level sets.

%
%
%
\begin{theorem}~\label{th16b} {\em  (Infinite Number of Local Level Sets)}
There exists a  dense set of ordinate values $y$  in $[0, \frac{2}{3}]$, which
are all dyadic rationals, 
such that  the level set $L(y)$ contains an infinite number of distinct local level sets. 
\end{theorem}

This theorem follows directly from a result proved  in Sect. \ref{sec5B} (Theorem \ref{th62a}). 
This in turn is derived  from the fact  that  $L(\frac{1}{2})$ is a countable set which 
contains a countably infinite number of  local level sets (Theorem \ref{th38}). \medskip

In the final Sect. \ref{sec9} we formulate some open questions
suggested by this work.


\subsection{Extensions of results and related work} 

The Takagi function has self-affine properties,
and there  has been extensive study of 
  various classes of self-affine functions. 
In particular, in the late 1980's Bertoin
\cite{Ber88}, \cite{Ber90} studied
 the Hausdorff dimension
of  level sets of certain classes of self-affine functions;
however his results do not cover the Takagi function.

In 1997 Yamaguti, Hata and Kigami \cite[Chap. 3]{YHK97} gave
a general definition of  a family $F(t, x)$ of Takagi-like functions
 depending on a parameter $0 < t < 1$ as follows:
 let $g(x)$ be a bounded measurable function defined on $[0,1]$
and let $\Phi: [0, 1] \to [0,1]$ be a continuous mapping, then set 
\begin{equation*}
F(t, x) := \sum_{n=0}^{\infty} t^n g( \Phi^n(x)).
\end{equation*}
If we specialize to take $\Phi(x)= 2 g(x)= T(x)$, the tent map in (\ref{eq101a}),
then the parameter value $t=\frac{1}{2}$ gives the Takagi function
$$
F(\frac{1}{2}, x) = 2 \left(\sum_{n=1}^{\infty} \frac{1}{2^n} T^{n}(x) \right)= \tau(x).
$$
If one now changes the parameter value to $t= \frac{1}{4}$,
then one gets instead (\cite[p. 35]{YHK97}) the smooth  function
\begin{equation*}
F(\frac{1}{4}, x) = \frac{1}{2}\left( \sum_{n=1}^{\infty} \frac{1}{4^n} T^{n}(x) \right)= \frac{1}{2}x(1-x).
\end{equation*}
The level sets of this function are finite. These examples  show
an extreme dependence of level set structure on the parameter $t$.
 Our analysis 
uses the piecewise linear nature of the function $\Phi(x)$
in a strong way, and also uses specific properties of the geometric scaling 
by the parameter $t$ at value $t=\frac{1}{2}$. \smallskip

The methods presented should extend to various functions  similar in construction to
the Takagi function, such as van der Waerden's function (\cite{vdW30}). 
They also  extend to intersections
of the graph of the Takagi function with parallel families of lines having
integer slope, a device used by the second author \cite{Maddock09}. 
In this paper we have treated only Hausdorff dimension, while the paper
 \cite{Maddock09}  also obtained upper bounds for Minkowski 
 dimension (a.k.a. box counting dimension) of level sets. Some results of this
 paper (e.g. Theorem \ref{th15}) may be strengthened to give
 Minkowski dimension upper bounds. \smallskip

In \cite{LM10b}
we further analyze the structure of global level sets $L(y)$ using local
level sets. We give a new proof of  a theorem of Buczolich \cite{Buz08} showing
that  if one draws   $y$ uniformly from
$[0, \frac{2}{3}]$, then with probability one the  level set $L(y)$  is a finite set; 
we  improve on it by
showing that the expected number of points in such  a ``random" level set $L(y)$
is infinite.  We also complement this result by  showing  that the
set of levels $y$ having a level set of positive Hausdorff dimension is
``large" in the sense that it has full Hausdorff dimension $1$,
although it is of Lebesgue measure $0$.\smallskip

Subsequent to this paper, Allaart \cite{A11} \cite{A11b} 
obtains many further results on local level sets. He shows that
dyadic ordinates $y=\frac{k}{2^n}$ have finite or countable level sets, and he
determines information on cardinalities of finite level sets.

Finally we remark that
there has been much study of the non-differentiable nature
of the Takagi function in various  directions, see for example
Allaart and Kawamura (\cite{AK06}, \cite{AK10})  and references therein.
It  is considered  as an example  in Tricot \cite[Section 6]{Tri97}. \smallskip

%
%
%
%
\section{Basic Properties of the Takagi Function}\label{sec2}
\setcounter{equation}{0}
We recall 
some basic facts and include proofs for the reader's convenience.
We first  give Takagi's formula for his function, which 
assigns a value $\tau(x)$ directly
to a  binary expansion of $x=0.b_1b_2b_3...$.
Dyadic rationals $\frac{k}{2^n}$ have two distinct binary expansions, and
one checks the assigned value $\tau(x)$
is the same for both expansions.
For $0 \le x \le 1$ the distance to the nearest integer function $\ll x \gg$ is
\begin{equation*}
\ll x \gg  ~:= \left\{ 
\begin{array}{lcl}
x & \mbox{if}& 0 \le x  < \frac{1}{2},~~\mbox{i.e.} ~b_1=0\\
~&~&~\\
1-x & \mbox{if}&  \frac{1}{2} \le x \le 1, ~~\mbox{i.e.} ~ b_{1}=1.
\end{array}
\right.
\end{equation*}
For  $n \ge 0$, we have 
\beql{204b}
\ll 2^n x \gg ~= \left\{ 
\begin{array}{lcl}
0.b_{n+1} b_{n+2} b_{n+3} ... & \mbox{if}& b_{n+1}=0\\
~&~&~\\
0.\tb_{n+1} \tb_{n+2} \tb_{n+3} ...& \mbox{if}&   b_{n+1}=1,
\end{array}
\right.
\eeq
where we use the bar-notation
\begin{equation*}
\bar{b}= 1-b , ~~~\mbox{for}~~~ b=0 ~\mbox{or} ~1,
\end{equation*}
to mean complementing a bit.

%
%
%
%

\begin{lemma}~\label{le21} {\rm (Takagi \cite{Takagi})}
For $x= 0.b_1b_2 b_3 ...$ the Takagi function is given by
\beql{201}
\tau(x) = \sum_{m=1}^{\infty} \frac{\ell_m}{2^m},
\eeq
in which $0 \le \ell_m= \ell_m(x) \le m-1$ is the integer
\begin{equation*}
\ell_m(x) = \# \{ i:~ 1 \le i < m,~~b_i \ne b_{m} \}.
\end{equation*}
In terms of the digit sum function $\dM_m(x)= b_1 + b_2 + ...+b_m$,
\beql{203}
\ell_{m+1} (x) = \left\{ \begin{array}{lcl}
\dM_{m}(x) & \mbox{if}& b_{m+1}=0,\\
~&~&~\\
m-\dM_{m}(x) & \mbox{if}& b_{m+1}=1.
\end{array}
\right.
\eeq
\end{lemma}

\paragraph{Proof.} From the definition 
\begin{equation*}
\tau(x) = \sum_{n=0}^{\infty} \frac{ \ll 2^n x \gg}{2^n}
\end{equation*}
Now \eqn{204b} gives
\begin{equation*}
\frac{\ll 2^n x \gg}{2^n} =
 \left\{ 
\begin{array}{lcl}
\sum_{j=1}^{\infty} \frac{b_{n+j}}{2^{n+j}} & \mbox{if}& b_{n+1}=0,\\
~&~&~\\
 \sum_{j=1}^{\infty} \frac{\tb_{n+j}}{2^{n+j}} & \mbox{if}&   b_{n+1}=1.
\end{array}
\right.
\end{equation*}
We substitute this into the formula for $\tau(x)$ and collect all
terms having a given denominator $\frac{1}{2^{m}}$, coming from $m=n+j$ with
$1 \le j \le m.$
For $m= n+j$ we get a contribution of $\frac{1}{2^m}$ whenever $b_{n+j} :=b_{m}= 1$ if $b_{n+1}=0$,
and whenever $b_{n+j}:=b_m=0$ if $b_{n+1}=1$, otherwise get $0$ contribution. Adding up over $j$,
we find  the total contribution is $\frac{\ell_m}{2^m}$ where $\ell_m(x)$ counts the number of $b_j$,
$1 \le j < m$ having the opposite parity to $b_m$, which is \eqn{201}. 
Note that $\ell_1(x) \equiv 0$, so the summation \eqn{201} really starts with $m=2$.
The formulas \eqn{203} follow by inspection; note that 
$m- \dM_m(x) = \dM_m(1-x)$ (making an appropriate convention for dyadic rationals).
$~~~\Box$\\

We next recall two basic functional equations,
see Kairies, Darslow and Frank  \cite{KDF88}.
%
%
%
%

\begin{lemma}~\label{le22}{\em (Takagi functional equations)} \\
The Takagi function satisfies  two functional equations,
each valid for 
$0 \le x \le 1.$ 
These are the reflection equation
\beql{206a}
\tau(x) = \tau(1-x),
\eeq
and the dyadic self-similarity equation
\beql{206b}
2 \tau(\frac{x}{2})= \tau(x) + x.
\eeq
\end{lemma}

\paragraph{Proof.} 
Here \eqn{206a} follows directly from   \eqn{eq101}, since
$\ll k x\gg ~= ~\ll k(1-x) \gg$ for $k \in \ZZ$.
To obtain \eqn{206b}, let $x= 0. b_1 b_2 b_3 ...$  and set
$y := \frac{x}{2} = 0.0 b_1 b_2b_3 ...$.  Then $\ll y\gg = y$, whence \eqn{eq101} gives
$$
2\tau(y) 
 = 2 \ll y \gg + 2\left(\sum_{n=1}^{\infty} \frac{ \ll 2^{n} y\gg} {2^n}\right) 
=x + \sum_{m=0}^{\infty}  \frac{ \ll 2^{m} x\gg} {2^m} = x + \tau(x).
~~~~\Box$$

We note  that the Takagi function 
can be characterized as the unique continuous
function on $[0,1]$ satisfying these two functional equations (Knuth \cite[Exercise 82, solution p. 740]{Kn09}).

Next we recall the construction of the Takagi function $\tau(x)$
as a limit of piecewise linear approximations
\begin{equation*}
\tau_n(x) := \sum_{j=0}^{n-1}  \frac{ \ll 2^j x\gg}{2^j},
\end{equation*}
which we name the {\em partial Takagi function} of level $n$.
We  require some notation concerning 
the binary expansion:
\begin{equation*}
x= \sum_{j=1}^{\infty} \frac{b_j}{2^j}= 0.b_1 b_2 b_3..., ~~~~ \mbox{each}~b_j \in \{0, 1\}. 
\end{equation*}
%
\begin{defi}\label{de23} 
{\em 
Let $x \in [0,1]$ have binary expansion
$x= \sum_{j=1}^{\infty} \frac{b_j}{2^j}= 0.b_1 b_2 b_3...$, with each $b_j \in \{0, 1\}$. 
For each $j \ge 1$ we define the following integer-valued functions.

(1) The {\em digit sum function} $\dM_j(x)$ is
\begin{equation*}
\dM_j(x) := b_1 +b_2 + \cdots + b_j.
\end{equation*}
 We also let $N_j^{0}(x)  = j - \dM_j(x)$ count the number of $0$'s in the 
 first $j$ binary digits of $x$.
 
 \smallskip
(2) The {\em  deficient digit function} $\dN_j(x)$ is given by
\beql{222}
\dN_j(x):= N_j^{0}(x) - \dM_j(x) = j- 2\dM_j(x)  = j-  2(b_1+b_2+ \cdots + b_j) .
\eeq
Here we use  the convention that $x$ denotes a binary expansion; 
dyadic rationals have two different binary expansions, and all functions 
$N_j^0(x)$, $\dM_j(x)$, $\dN_j(x)$ depend
on which binary expansion is used.
}
\end{defi}

The name ``deficient digit function" reflects the
fact that  $\dN_j(x)$ counts the excess of binary digits $b_k=0$ over those with $b_k=1$
in the first $j$ digits, i.e. it is positive if there are more $0$'s than $1$'s.

%
%
%
%

\begin{lemma}~\label{le23}{\em (Piecewise linear approximations to Takagi function)} \\
The piecewise linear function
$\tau_n(x) = \sum_{j=0}^{n-1}  \frac{ \ll 2^j x\gg}{2^j}$
 is  linear on each dyadic
interval $[\frac{k}{2^n}, \frac{k+1}{2^n}]$.

(1) On each such interval $\tau_n(x)$ has integer slope between $-n$ and $n$ given by
the deficient digit function  
$$
\dN_n(x) = N_{n}^0(x) - \dM_n(x) =n - 2(b_1+b_2 + \cdots + b_n),
$$
Here  $x=0.b_1b_2 b_3...$ may be any interior point on the 
dyadic interval, and can also be an endpoint provided the dyadic
expansion ending in $0$'s is taken at the left endpoint
$\frac{k}{2^n}$ and that ending in $1's$ is taken
at the right endpoint $\frac{k+1}{2^n}.$  

 (2) The values $\{ \tau_n(x):  n \ge 1\}$
converge uniformly to $\tau(x)$, with
\beql{253a}
|\tau_n(x) - \tau(x)| \le \frac{2}{3}\cdot\frac{1}{2^{n}}.
\eeq
The functions $\tau_n(x)$ approximate the
Takagi  function monotonically from below
\begin{equation*}
\tau_1(x) \le \tau_2(x) \le \tau_3(x) \le ...
\end{equation*}
For a dyadic rational $x= \frac{k}{2^n}$, perfect approximation
occurs at the $n$-th step with
\begin{equation*}
\tau(x) = \tau_m (x), ~~~\mbox{for all}~~ m \ge n.
\end{equation*}
\end{lemma}

\paragraph{Proof.} 
All statements follow easily from the observation that 
each function
$f_n(x) := \frac{ \ll 2^n x\gg}{2^n}$
is a piecewise linear sawtooth function, linear on dyadic intervals
$[\frac{k}{2^{n+1}}, \frac{k+1}{2^{n+1}}]$, with slope 
having value $+1$ if the binary expansion of $x$ has $b_{n+1}=0$
and slope having value $-1$ if $b_{n+1}=1$.  The inequality in
\eqref{253a} also uses the fact that $\max_{x \in [0,1]} \tau(x) =
\frac 2 3.$ 
$~~~\Box$\\

The Takagi function itself can be directly
expressed in terms of the deficient digit function.
The relation \eqn{222} compared with the definition \eqn{203} of 
$\ell_m(x)$ yields
$$
\ell_{m+1}(x) = \frac{m}{2}-  \frac{1}{2} (-1)^{b_{m+1}}\dN_{m}(x).
$$
Substituting this in Takagi's formula \eqn{201} 
and simplifying yields the formula
\begin{equation*}
\tau(x)
= \frac{1}{2} -\frac{1}{4} \left(\sum_{m=0}^{\infty} (-1)^{b_{m+1}} \frac{\dN_{m}(x)}{2^m}
\right).
\end{equation*}

We conclude this section with a self-similarity property of the Takagi
function associated to dyadic rationals $x = \frac{k}{2^{n}}$.

%
%

\begin{lemma}~\label{le25}{\em (Takagi self-affinity)} \\
For an arbitrary dyadic rational $x_0=\frac{k}{2^n} $ then for  $x \in [\frac{k}{2^n}, \frac{k+1}{2^n}]$
given by $x = x_0 + \frac{w}{2^{n}}$ 
with $w\in [0,1]$,
\begin{equation}
\label{251}
\tau( x) = \tau(x_0) + \frac{1}{2^n} \large( \tau(w) + \dN_n(x_0)w\large).
\end{equation}
That is, the graph of $\tau(x)$ on  $ [\frac{k}{2^{n}}, \frac{k+1}{2^{n}}]$
is a miniature version of the tilted Takagi function
$ \tau(x) + D_n(x_0) x$, 
vertically shifted by $\tau(x_0)$, and  shrunk by a factor
$\frac{1}{2^{n}}$. 
\end{lemma}

\paragraph{Proof.} 
By Lemma \ref{le23}(1), we have $\tau_{n}(x_0 + \frac{w}{2^{n}}) =
\tau_{n}(x_0) + \dN_n(x_0)\cdot \frac{w}{2^n}$. 
Therefore, by \eqref{eq101} it follows that
\begin{align*}
  \tau(x) & = \tau_{n}(x) + \sum_{j=n}^\infty
  \frac{\ll 2^j x\gg}{2^k}\\
  & = \tau_{n}(x_0) + \dN_n(x_0)\cdot \frac{w}{2^n} +
    \sum_{j=n}^\infty \frac{\ll 2^j (\frac{w}{2^{n}}) \gg}{2^j}\\
  & = \tau(x_0) + \frac{1}{2^n} \large( \tau(w) + \dN_n(x_0) w \large).~~~\Box.
\end{align*}

\paragraph{Remark.}
 Lemma \ref{le25} simplifies 
 in the special case that $D_n(x_0)=0$, in which case
  we call $x_0$ a  {\em  balanced dyadic rational};
 this can only occur when  $n=2m$ is  even.  The formula \eqn{251} becomes
\begin{equation*}
\tau( x) = \tau(x_0) + \frac{\tau(w)}{2^n},
\end{equation*}
which shows the graph of the Takagi function over the subinterval $[\frac k
  {2^n}, \frac{k+1}{2^n}] \subseteq [0,1]$, up to a translation,
consists of the image of the entire graph scaled by $\frac 1 {2^n}$.
Balanced dyadic rationals play a special role in our analysis,
see Definition \ref{de32}.

%
%
%
%
\section{Properties of Local Level Sets}\label{sec3}
\setcounter{equation}{0}

In this section we derive some basic properties of local level sets.
Then we determine the size of ``abscissa generic" local level sets,
showing these are uncountable sets of Hausdorff dimension $0$.
Finally we determine the structure of local level sets that contain
a rational number $x$, showing they are either finite sets or 
Cantor sets of positive Hausdorff dimension.


\subsection{Partition into local level sets} \label{sec31}

We first show that level sets partition into local level sets.
%
%

\begin{theorem}~\label{th31}
(1) Local level sets  $L_x^{\loc}$ are closed sets. Two local level sets
either coincide or are  disjoint. 

(2) Each local level set $L_x^{\loc}$ is contained in a level set:
 $L_x^{loc} \subseteq L(\tau(x))$.
That is, if  $x_1 \sim x_2$ then $\tau(x_1) = \tau(x_2).$

(3) Each level set $L(y)$ partitions into local level sets
\begin{equation}
\label{300e}
L(y) = \bigcup_{{x \in \ddL}\atop{\tau(x)= y}}  L_{x}^{\loc}.
\end{equation}
Here $\ddL$ denotes the collection of leftmost endpoints of all local level sets.
\end{theorem}

\paragraph{Proof.} 
(1) A  local level set, specified by a binary expansion of $x$, is 
 generated by any one of its elements,
by  block-flipping operations (allowing an infinite
number of blocks to be flipped at once) including
$x \to 1-x$ : if $x_2 \in L_{x_1}^{\loc}$
then $x_1 \in L_{x_2}^{\loc}$ and $L_{x_1}^{loc}= L_{x_2}^{\loc}$. 
Each set $L_x^{\loc}$ is closed, 
because any Cauchy sequence $\{ x_n: n \ge 1\}$
in $L_x^{\loc}$ must eventually ``freeze" the choice made in any finite
initial set of  ``blocks", so that for each $k \ge 1$ there is a value
$n(k)$ so that
$B_k(x_n)= B_k(x_m)$ whenever $n,m \ge n(k)$. Since 
all balance-sets  $Z(x_n)$ coincide, the limit value $x_{\infty}$
has $B_k(x_{\infty}) = B_k(x_n)$ for all $n \ge n(k)$ so
$B_k(x_{\infty}) \sim B_k(x)$ for all $k \ge 0$. This same relation shows
that the first $k$ balance points of $x_{\infty}$ coincide with those of such $x_n$,
and letting $k \to \infty$ we have $Z(x_{\infty}) = Z(x)$. Thus $x_{\infty} \sim x$,
so $x_{\infty} \in L_x^{\loc}$. \\

(2) We assert that if $x \sim x'$, with $x'$ obtained from $x$ by flipping a single block of symbols,
then $\tau(x) = \tau(x')$.  This holds since
a block-flip after the $k$th binary digit of $x$ corresponds to
a reflection of $x$ about the center of the dyadic interval of length
$\frac 1 {2^k}$ containing $x$;  by Lemma \ref{le25} and Lemma \ref{le22}, the Takagi
function restricted to this interval has this reflection symmetry.
   The case $\tau(x) = \tau(x'')$ of general $x''$ 
in $L_x^{loc}$ will then follow by 
flipping
the blocks of $x$ in increasing order as necessary
to match those of $x''$, getting a sequence $\{ x_n: n \ge 1\} \subset L_x^{\loc} $
with $\tau(x_n) = \tau(x)$. Now $\lim_{n \to \infty} x_n = x''$ 
so  using the fact that $\tau(x)$ is a continuous
function, we conclude $\tau(x'') = \lim_{n \to \infty} \tau (x_n) = \tau(x).$\\

(3)   Local  level sets, being closed,  have a leftmost
 endpoint, and we 
can then uniquely label local level sets with their leftmost endpoint. $~~~\Box$ \medskip

We immediately deduce that any level set that is countably infinite must
contain infinitely many local level sets. 
%
%

\begin{coro}~\label{cor32}
 {\em(Countable Level Sets)}
 (1) Each local level set $L_x^{\loc}$  is either  finite or
 uncountable.

(2) Each level set $L_x^{\loc}$ that is countably infinite
is necessarily a countable disjoint union of finite local level sets.
\end{coro}

\paragraph{Proof.} 
(1) This dichotomy for $L_x^{loc}$ is determined by whether the balance set $Z(x)$ is
finite or infinite, since $L_x^{loc}$ is a Cantor set in the latter case.

(2) If $L(y)$ is countably infinite, all local level sets it contains must
be finite, by (1). Thus there must be infinitely many of them.
$~~~\Box$\\

We will later show that
case (2) above occurs:  Theorem \ref{th38} proves  that 
$L(\frac{1}{2})$ is countably infinite.


\subsection{Generic local level sets}\label{sec32}

We analyze the size of  ``abscissa generic" local level sets sampled by choosing
$x$ uniformly on $[0, 1]$, and prove these are 
Cantor sets of Hausdorff dimension $0$.

%
%

\paragraph{Proof of Theorem \ref{th15}.}
The sequence of binary digits in $x=0.b_1b_2b_3...$ of a real number
in $[0,1]$ corresponds to taking a  (random) walk on 
the integer lattice, starting at point $s_0=0$, with steps $+1$ or $-1$, 
with $b_j=0$ corresponding to taking a step in the positive direction,
and $b_j=1$ corresponding to taking a step in the negative direction, i.e.
at time $k$ the walk is at position
$$
s_k = s_0 + \sum_{j=1}^k (-1)^{b_j}.
$$
The interval $[0,1]$ sampled by drawing a random
 point  $x$ with the uniform distribution  (i.e. using Lebesgue measure on $[0,1]$)
  corresponds  in probability to taking a simple random walk with equal probability steps.
(See Billingsley \cite[Sect. 3]{Bil60},\cite{Bil65}).\\

The relevant property of the random walk is that 
a one-dimensional random walk is {\em recurrent}; that is, with probability one it
returns to the origin infinitely many times.
Thus with probability one the balance-set $Z(x)$ includes infinitely many
balance points, so that  with probability one the local level set $L_x^{loc}$ has the
structure of  a Cantor set, hence is uncountable. This corresponds to a
full Lebesgue measure set of points $x$ having a local level set that
is uncountable.\\

For a treatment of Hausdorff dimension, see Falconer  \cite{Fa03}.
To establish  the Hausdorff dimension $0$ assertion, we use the result
that with probability one  the number of times a simple random walk 
returns to the origin in the first $n$ steps is $o(n)$ as $n \to \infty$; in fact
with probability one  it is $O\left( n^{\frac{1}{2} + \epsilon}\right)$ as $n \to \infty$.
(See Feller \cite[Chap III]{Fe68}, \cite[Chap XII]{Fe71}).
The proof is then completed by the following deterministic result.

\paragraph{Claim.}
{\em Any $x=0.b_1b_2b_3...$  that has the
 property that  the number of returns to the origin in the first $n$ steps
of the corresponding random walk  is 
 $o(n)$ 
 as $n \to \infty$ necessarily has Hausdorff dimension
$\dim_{H} ( L_x^{\loc}) = 0.$ }\\

To show this, let the balance points of $Z(x)$ be
$0 < c_1< c_2 < c_3< ...$.  The hypothesis on $x$ implies 
$$
\lim_{ k \to \infty} \frac{k}{c_k} = 0.
$$
We can now cover  $L_x^{\loc}$ by 
 $2^k$-dyadic intervals of length $2^{-c_k}$, since there are only $2^k$ possible
 flipped sequences; call this covering $\sC_k$.  For any $\delta>0$, we have
 $$
 \lim_{k \to \infty} \sum_{I_j \in \sC_k}  |I_j|^\delta= 2^k 2^{- \delta c_k} = 0,
 $$
 since $\frac{c_k}{k} \to \infty.$ This proves the claim, 
 which completes the proof. 
  $~~~\Box$\\


\subsection{ Local level sets containing rational numbers} \label{sec33}

Knuth~\cite[Sect. 7.2.1.3, Exercise 83]{Kn09} raised the question of determining which
rational $y$ have an uncountable level set $L(y)$. 
We address here the easier question
 of  determining  which rational numbers $x$ 
 have an uncountable local level set $L_x^{loc}$.
We also  show that  uncountable local level sets  that contain a rational  
 necessarily have positive Hausdorff dimension.

%
%

\begin{theorem}~\label{th50} {\em (Rational local level sets)}

For a rational number $x=\frac{p}{q} \in [0,1]$, the following properties are equivalent.
\begin{enumerate}
 \item[(1)] The local level set $L_x^{loc}$ has positive Hausdorff dimension.
 
 \item[(2)] The local level set $L_x^{loc}$ is uncountable.
 
 \item[(3)] The  binary expansion of $x$  has a purely periodic part
 with an equal number of zeros and ones, and also has a
preperiodic part with an equal number of zeros and ones.
\end{enumerate}
Moreover, if these equivalent properties hold, then $\dim_H(L_x^{loc})
= \frac k r$, where $r$ is the number of bits in the periodic part of the binary
expansion of $x$ and $k$ is the number of balance points per period.
\end{theorem}

\paragraph{Proof.} 

Trivially, (1) implies (2).\medskip

To show  (2) implies (3), let $x = 0.b_1b_2b_3\ldots$ be the
binary expansion of $x = \frac{p}{q}$.  
Let  the balance-set $Z(x) := \{ c_j : \dN_{c_j}(x) = 0 \}$ be the set of balance points of $x$, 
as defined in
(\ref{123}).  By inspection, condition (3) is equivalent to the set $Z(x)$
having 
infinite cardinality.  From the definition of local level set 
$L_x^{loc}$, we see that the cardinality of $L_x^{loc}$ is 
$$ \# L_x^{loc} = 2^{\#Z(x)}.$$
Hence, $L_x^{loc}$ is uncountable if and only if (3) holds. \medskip

It remains to show that (3) implies (1).  
Since Hausdorff dimension is invariant under scaling, translation, and
finite union, we may reduce to the case where the rational $x$ has no
preperiodic part; that is, $x = 0.(b_1b_2\ldots b_r)^\infty$, with $r$ the
length of the periodic part.  Say
that these $r$-bits are partitioned into precisely $1 \leq k \leq
\frac r 2$
blocks by the balance points of $x$.
Let $B_1,\ldots, B_{2^k}$ be the $2^k$ possible dyadic rationals
obtained by applying block-flipping operations to $B_1 = 0.b_1\ldots
b_r$.
Define the function $S_i(x) := B_i + \frac{x}{2^r}.$  Using the
terminology of Falconer \cite[\S 9.2]{Fa03},
the functions $\{S_1, \ldots, S_{2^k} \}$ form an iterated function
system with attractor equal to the local level set 
$L_x^{\loc} = \bigcup_{i=1}^{2^k} S_i(L_x^{\loc}).$
Furthermore, the open set $V = (0,1)$ is a bounded open set for which 
$V \supseteq \bigcup_{i=1}^{2^k} S_i(V)$ with the union disjoint.
Therefore, it satisfies the hypothesis of  
\cite[Theorem 9.3]{Fa03}, whose conclusion yields that the Hausdorff
dimension 
$\dim_H(L_x^{\loc}) = s$
where $s$ is the unique number for which
$  \sum_{i=1}^{2^k}c_i^s = 1$,
and $c_i = 2^{-r}$ is the ratio of similitude of the operator $S_i$.
The equation is easily solved for $s$ yielding $s = \frac k r > 0$.
Thus, (3) implies (1).
(It is also a consequence of  \cite[Theorem 9.3]{Fa03} that the
$s$-dimensional Hausdorff measure of $L_x^{\loc}$ is a finite positive
value.) 
$~~~\Box$\\

\paragraph{Remarks.}
(1)  By Theorem \ref{th50}(3), any 
local level set of a rational number $x$ 
that is uncountable  necessarily 
contains infinitely many rational numbers $x'$,
by using periodic sequences of flippings.

(2) A  dyadic rational $x$ always belongs to a
finite local level set $L_x^{loc}$.  This follows immediately from
Theorem \ref{th50} since local level sets are either finite or
uncountable.

%
%
%
%
\section{Left Endpoints of Local Level Sets}\label{sec4}
\setcounter{equation}{0}

We study the set of endpoints of local level sets. The leftmost endpoints
fall in $[0, \frac{1}{3}]$ and the rightmost endpoints in $[\frac{2}{3}, 1]$,
and these sets are related by the operation $x \to 1-x$.
Therefore  suffices to study
the leftmost endpoint set, denoted $\ddL$. 
The usefulness of this set is that it  parametrizes
the complete collection of all local level sets,
and the Takagi function turns out to be well-behaved when
restricted to $\ddL$.
Theorem \ref{th33} 
below describes the main properties of $\ddL$.


\subsection{Deficient  digit  set} \label{sec41}

We make the following definition, which  will  be shown later to 
 coincide with the set of  all leftmost endpoints of local level sets.

%
%

\begin{defi}\label{de31}
{\em
The {\em deficient  digit  set} $\ddL$ consists of all
points
$$
\ddL := \left\{ x = \sum_{j=1}^{\infty} \frac{b_j}{2^j}:  ~~ \dN_j(x)
\ge 0 ~\mbox{for ~all} ~~j \ge 1\right\},
$$
in which  the deficient digit function $\dN_j(x)= j-2\dM_j(x)$ counts the number of binary
digits $0$'s minus the number of  $1$'s in the first $n$ digits. \smallskip
Note   that dyadic rationals $\frac{m}{2^r}$
have two different binary expansions; at most one of
these two expansions can belong to $\ddL$;  
if one of them does, then by convention
we assign the dyadic rational to the set $\ddL$.
}
\end{defi}

We will establish in Theorem~\ref{th33}  that the deficient digit set $\ddL$ is a Cantor set
having Lebesgue measure zero. 
We need the following result as  a preliminary step.
Recall that a  dyadic rational binary expansion
$x:= \frac{k}{2^n}=0.b_1 b_2... b_n0^{\infty}$ (with $k$ odd)  is
said to be  {\em balanced} if
$\dN_{n}(x)=0$, and that  $n=2m$ is necessarily even.

%
%

\begin{lemma}\label{le31} {\em (Balanced dyadic rationals in $\ddL$)}
For a fixed integer $m\ge 0$,
the set  $\sBp$ of dyadic rationals $\frac{k}{2^{2m}}= 0.b_1 b_2 \cdots b_{2m}$ that 
are balanced and belong to $\ddL$, i.e. have  digit sums satisfying 
\begin{equation*}
\dN_{j}\left(\frac{k}{2^{2m}}\right)  \ge  0 ~~\mbox{for}~~1\le j \le 2m-1~~~\mbox{and}~~ 
\dN_{2m}\left(\frac{k}{2^{2m}} \right) =0,
\end{equation*}
has cardinality the  $m$-th Catalan number $C_{m} = \frac{1}{m+1} \left( {{2m}\atop{m}} \right),$
with  $C_0=1$. 
\end{lemma}

\paragraph{Proof.}
Each such dyadic rational describes  a lattice path starting from $(0,0)$ and
taking steps $(1,1)$ or $(1, -1)$ in such a way as to stay on or above the line $y=0$.
Here the $j$-th step $(1,1)$ corresponds to $b_j=0$ and $(1,-1)$ to $b_j=1$;
the last step necessarily has $b_{2m}=1$.
Such steps can be  counted using Bertrand's ballot theorem
(Feller \cite[p.73] {Fe68}).  To apply the theorem
we count  paths from $(0,0)$ to $(n, x) = (2m+1, 1)$ that 
stay  strictly above the $x$-axis.   Note that all such paths must go
through $(1,1)$, and 
that to end at $(2m+1, 1)$ the last step must be $(1, -1)$.  The
number of paths is therefore
$$
\frac{1}{2m+1} \left( {{2m+1}\atop{m}} \right)= \frac{1}{m+1}  \left( {{2m}\atop{m}} \right)=C_m,
$$
as asserted.
$~~~\Box$\\

We will determine the set of open intervals removed from $[0,1]$
to create the deficient digit sum set $\ddL$. The following set $\sB$ supplies
labels for these open intervals.

%
\begin{defi} \label{de32}
{\em

(1) The {\em breakpoint set}
 $\sBp$ is the set of
 all balanced dyadic rationals in $\ddL$.
 It consists  of $\sBo'=0$ together with the  collection of  all dyadic rationals 
   $B' = \frac{n}{2^{2m}}$  that have  binary expansions of the form
$$
B'=0.b_1 b_2 ... b_{2m-1}b_{2m} ~~~\mbox{for some}~~ m \ge 1,
$$
that satisfy  the condition
\begin{equation*}
\dN_j(B') \ge 0 ~~~\mbox{for} ~~1 \le j \le 2m-1,
~~~\mbox{and}~~~\dN_{2m}(B') =0. 
\end{equation*}

(2) The {\em small breakpoint set} $\sB$ is the subset of 
the breakpoint set $\sB'$ consisting of 
$\sBo=0$ plus all members of $\sB'$ satisfying the extra
condition that  
the last two binary digits $b_{2m-1}= b_{2m}=1.$

\noindent We may rewrite  a dyadic rational in the small breakpoint set as
\begin{equation}
\label{332a}
B= 0. b_1 b_2 ... b_\ell 0 1^k, ~~~\mbox{with}~~ k \ge 2,
\end{equation}
where  $2m= k+\ell+1.$
}

\end{defi}

We show that values in the small breakpoint set $\sB$ 
naturally label  the left endpoints $x(B)^{-}$
of the intervals removed from
$[0,1]$ to create the deficient digit set $\ddL$.
(It is also possible to give a nice labeling for the right endpoints $x(B)^{+}$
which we omit here.)

%

\begin{defi} \label{de33b}
{\em
For each dyadic rational  $B=0.b_1b_2 ...b_{\ell} 01^k$, $k \ge 2$ 
 in the small breakpoint set $\sB$ ($B \ne \sBo$) we associate
 the open interval
\begin{equation*}
I_B := (x( B)^{-}, x(B)^{+})
\end{equation*}
 having the endpoints
\begin{eqnarray*}
x(B)^{-} & := &0. b_1 b_2 ... b_\ell 0 1^k (01)^{\infty}\\
x(B)^{+} & := & 0. b_1 b_2 ... b_{\ell} 1 0^k (00)^{\infty}.
\end{eqnarray*}
For $B= \sBo$ we set
\begin{equation*}
I_{\sBo} := ( x(\sBo)^{-} , x(\sBo)^{+}) :=( 0.(01)^{\infty},
1.(00)^{\infty}) = \left( \frac{1}{3}, 1\right). 
\end{equation*}
} 
\end{defi}

Some data on $I_{\sB}$ for small $\ell$ and $k$ appears in
Table \ref{tab2}. 

%
%
\begin{table}[h]
\begin{center}
$$
     \begin{array}{|c|c|c|c|c|c|}
\hline
~& B = j/2^{2m} & x(B)^{-} & x(B)^{+}  &  \tau(x(B)^{-})     & \tau(x(B)^+)\\
\hline
2m=4&&&&&\\
\hline
&3/16=.0011& 5/24=.0011(01)^\infty& 1/4= .0100&  13/24  &1/2\\
\hline
2m=6&&&&&\\
\hline
&7/64= .000111& 11/96= .000111(01)^\infty& 1/8 = .001000& 37/96 &3/8\\
\hline
&11/64=.001011& 17/96= .001011(01)^\infty& 3/16= .001100& 49/96&1/2\\
\hline
&19/64=.010011&29/96= .010011(01)^\infty& 5/16= .010100&  61/96&5/8\\
\hline
\end{array}
$$

  \end{center}
  \caption{Binary expansions of $B$, $x(B)^{-},$ and
    $x(B)^{+}$ for $B \in \sB$ of the form $B = \frac{j}{2^{2m}}$ 
    and corresponding $\tau(x(B)^{-})$ and $\tau(x(B)^{+})$. } 
      \label{tab2}
\end{table}

%
%

\begin{lemma}~\label{le32}
The open intervals $\{ I_B:~ B \in \sB\}$ have the following properties.\\

(1) The intervals  $I_B$ for $B \in \sB$ are all disjoint from the deficient digit set $\ddL$
and from each other. All the endpoints $x(B)^{\pm}$ 
belong to $\ddL$, with the exception of 
$x(\sBo)^{+}= 1$. \\

(2) For $B$ in the small breakpoint set $\sB$  there  holds 
\begin{equation*}
x(B)^{+} - x(B)^{-}  = \tau(x(B)^{-})- \tau(x(B)^{+})= \frac{1}{2^{k+\ell} \cdot 3} =\frac{1}{2^{2m-1} \cdot 3}
\end{equation*}
so that $x(B)^{+} > x(B)^{-}$. Thus the ratio $\frac{\tau(x(B)^{+})-
  \tau(x(B)^{-})} {x(B)^{+} - x(B)^{-}} = -1$.
\end{lemma}

\paragraph{Proof.}
(1) We have  $2m= k+\ell +1$ in \eqn{332a}, and with the 
definition of $x(B)^{+}, x(B)^{-}$ this implies that there
are odd integers $n_1, n_2$ such that
\begin{eqnarray*}
x(B)^{-} &= &\frac{n_1}{2^{k+\ell+1}} + \frac{1}{2^{k+\ell+1} \cdot 3}\\
x(B)^{+} &=& \frac{n_2}{2^{\ell+1}}.
\end{eqnarray*}
Furthermore one easily sees that $x(B)^{+} > x(B)^{-}$ and that
\begin{equation}
\label{338c}
x(B)^{+} - x(B)^{-} = \frac{ 2}{2^{k+\ell+1}\cdot 3} = \frac{1}{2^{2m-1}\cdot 3}.
\end{equation}
A key property of $x(B)^{-}$ is that it belongs to $\ddL$ and satisfies
$$
\dN_{k+\ell+1+ 2j} (x(B)^{-}) = 0 ~~~\mbox{for all} ~~j \ge 0.
$$
Now the binary expansion of any number $x$ strictly between $x(B)^{-}$
and $x(B)^{+}$ 
first differs from $x(B)^{-}$ in a bit at location $\ell' > \ell$,
with $x$ having digit $1$ and $x(B)^{-}$ digit $0$. But a $1$ in digit
location   
$k+\ell + 2j, ~~j \ge 1$  would then produce
$$
\dN_{k+\ell+ 2j} (x)= -1,
$$
which certifies $x \not\in \ddL$. If  there is instead
a change from $0$ to $1$ in a digit of $x(B)^{-}$  in location $j \le \ell+1$,
then this would make $x \ge x(B)^{+}$, a contradiction. We conclude that
no point of the open interval $I_B$  belongs to   $\ddL$.

We note by inspection 
that all upper endpoints $x(B)^{+} \in \ddL$, except $x(\sBo)^{+}=1$.
This certifies that all open intervals $I_B$
are disjoint (except possibly $I_{B_{\emptyset}}$) , because the endpoints of the closure of each 
such interval are in $\ddL$ but the interiors are not. 
(This  prevents
both overlap and inclusion.) Finally, $I_B \subset (0, \frac{1}{3})$
for each interval with $B \ne \sBo$, yielding disjointness in all cases. \\

(2) In view of \eqn{338c} it remains to verify that
\beql{339a}
\tau(x(B)^{-}) - \tau( x(B)^{+} ) = \frac{1}{2^{k+\ell}\cdot 3}.
\eeq
In the notation of Lemma~\ref{le21} we write
$$
y(B)^{+} :=\tau( x(B)^{+} ) = \sum_{j=1}^{\infty} \frac{\ell_j^{+}}{2^j}; ~~~~~~~
y(B)^{-} := \tau(x(B)^{-}) = \sum_{j=1}^{\infty} \frac{\ell_j^{-}}{2^j}.
 $$
 Now $\ell_j^{-} = \ell_j^{+}$ for $1 \le j \le \ell$, since the two expansions agree.
 We also have $\ell= 2m-k-1$ and find that for $B\in\sB$
 the  first $\ell$ digits of $B$ necessarily contain 
 $m_1 = m-k = \frac{\ell-k+1}{2}$ values $b_j=1$
 and $\frac{\ell+k-1}{2}$ values $b_j=0$. Now we calculate using Lemma~\ref{le21} that
 $$
 \ell_{\ell+1}^{-}= \frac{\ell-k+1}{2},  ~\mbox{and}~~ \ell_{\ell+1}^{+}= \frac{\ell+k-1}{2}
 $$
 while for $2 \le j \le k+1$ we have
 $$
\ell_{\ell+j}^{-} = \frac{\ell+k+1}{2},~~\mbox{and} ~~~\ell_{\ell+j}^{+} = \frac{\ell-k+3}{2}.
 $$
 Thus in the first $2m=k+\ell+1$ terms we have
 \beql{339c}
 \sum_{j=1}^{k+\ell+1} \frac{\ell_j^{-}- \ell_j^{+}}{2^j} = -\left(\frac{k-1}{2^{\ell+1}}\right)+ \sum_{j=1}^k \frac{k-1}{2^{\ell+1+j}}
 = \frac{k-1}{2^{\ell}} \left( -\frac{1}{2} + \sum_{j=2}^{k+1} \frac{1}{2^j} \right) = -\frac{ k-1}{2^{k+\ell+1}}.
 \eeq
For the terms $j \ge k+\ell+1$, for $x(B)^{+}$ we have $(m-k+1)$ $1$'s in the first $2m=k+\ell+1$
positions, hence
\beql{339d}
\sum_{j=2m+1}^{\infty} \frac{\ell_j^{+}}{2^j}= \sum_{j=2m+1}^{\infty} \frac{m-k+1}{2^j} = 
\frac{m-k+1}{2^{k+\ell+1}}.
\eeq
For  $x(B)^{-}$, we have equal numbers of $0$'s and $1's$ at the $2m$-th digit,
so that $\ell_j^{-} = 2m + \tilde{\ell_j}$ where $\tilde{\ell_j}$ correspond to the expansion
of $\frac{1}{3} = 0. (01)^{\infty}$. Since we have $\tau(\frac{1}{3}) = \frac{2}{3}$,
we obtain
\beql{339e}
\sum_{j=2m+1}^{\infty} \frac{\ell_j^{-}}{2^j}= \sum_{j=2m+1}^{\infty} \frac{m}{2^j} +\sum_{j=2m+1}^{\infty}
\frac{ \tilde{\ell}_j}{2^{j}} = \frac{m}{2^{2m}} + \frac{2}{3} \cdot \frac{1}{2^{2m}}
= \frac{ m+ \frac{2}{3}}{2^{\ell+k+1}}.
\eeq
Combining \eqn{339c}-\eqn{339e} yields
$$
\tau(x(B)^{-} ) - \tau(x(B)^{+}) = \frac{1}{2^{k+\ell+1}} \left( -(k-1) +(m+\frac{2}{3}) -(m-k+1)\right)
= \frac{2}{3} \cdot \frac{1}{2^{\ell+k+1}}= \frac{1}{2^{k+\ell} \cdot 3},
$$
verifying \eqn{339a}.
$~~~\Box$\\


\subsection{Properties of  the deficient digit  set} \label{sec42}

The following result characterizes the deficient digit set $\ddL$.

%
%

\begin{theorem}\label{th33} {\em(Properties of the deficient digit set)}\\
$~~~~~$(1) The deficient digit  set $\ddL$ comprises the set of leftmost endpoints of all local
level sets. It satisfies $\ddL \subset [0, \frac{1}{3}]$.

\smallskip
(2) The deficient digit set $\ddL$ is a closed, perfect set (Cantor set).
It is given by 
\beql{342}
\ddL= [0, 1) \backslash \bigcup_{ B \in \sB} I_{B},
\eeq
where the omitted open intervals $I_B$ have  right endpoint  a dyadic rational and 
left endpoint a rational number with denominator $3 \cdot 2^k$ for some $k \ge 1$.

\smallskip
(3) The deficient digit set $\ddL$ has Lebesgue measure zero.
\end{theorem}

\paragraph{Proof.}
(1) This property is immediate from the definition of local level
set. The leftmost 
endpoint of any local level set satisfies $\dN_j(x) \ge 0$ for all $j \ge 1$ and is the only point in
$L_x^{loc}$ with this property. \medskip

(2) The definition of $\ddL$ shows that it is a closed set, since
the inequalities $\dN_j(x) \ge 0$ are preserved under pointwise limits.
Note here that any infinite binary expansion ending $.11111...$ is excluded from
membership on $\ddL$. To see that $\ddL$ is a perfect set, we show each member of $\ddL$ is a limit
point of other members of $\ddL$.
For each member $x \in \ddL$ whose binary expansion contains an infinite number of $1$'s,
we approximate it from below by
the sequence $x_n= \frac{1}{2^n} \lfloor 2^n x\rfloor \in \ddL$ obtained by truncating it at the $n$-th digit. 
For a dyadic
rational member $x= \frac{k}{2^j}$, which necessarily ends in an infinite string of zeros, 
we  approximate it from above using the sequence $x_n = x+ \frac{1}{2^{n+j + 2}} \in \ddL$. 

To show the equality  \eqn{342} set
$\ddL_C := [0, 1) \backslash \bigcup_{ B \in \sB} I_{B}.$
We clearly have $\ddL \subseteq \ddL_C$, by Lemma~\ref{le32}(1).
It remains to show $\ddL_ C \subseteq \ddL$. We check the
contrapositive, that $x \not\in \ddL $ implies $x \not\in \ddL_C$. We use the criterion  
that if $x \not\in \ddL$
then $\dN_j(x) <0$ for some $j\ge 1$. Now one can verify that
the removed intervals $I_B$
each detect those $x$ whose first occurrence of
$\dN_j(x) <0$ is in a specified digit position $j$, with a specified digit pattern 
of the first $k$ digits, followed by some string of $(01)^r$, and these enumerate all
possibilities of this kind. Thus $x \not\in \ddL_C$, showing \eqn{342}.
The properties of the endpoints of $I_{B}$ are given in Lemma  \ref{le32}(2).
\medskip

(3) The set $\ddL$  is shown to have Lebesgue measure $0$ by
covering it with dyadic boxes at level $2m$ each of size $\frac{1}{2^{2m}}$,
and noting from Lemma~\ref{le31} that exactly $C_m$ such boxes need to be
used to cover $\ddL$, so that
$$
\meas(\ddL) \le \frac{C_m}{2^{2m}}.
$$
Stirling's formula gives for the Catalan numbers $C_m$ the estimate
\begin{equation*}
C_m =\left( 1 + o(1)\right) \frac{1}{\sqrt{\pi m^3} }~2^{2m} , ~~\mbox{as}~~ m \to \infty.
\end{equation*}
From this  we see that 
$ \frac{C_m}{2^{2m}} \to 0$ as $m \to \infty.$
$~~~\Box$.\\

%
%
%

\subsection{Takagi function on  deficient digit set} \label{sec44}

We next consider  the Takagi function restricted to the deficient
digit set $\ddL$. We show firstly that it has a weak
increasing property approaching any point of $\ddL$, and secondly that it is nondecreasing    
when further restricted to the   set 
 $\frac{1}{2} \ddL:= \{ \frac{1}{2} x:~ x \in \ddL\}.$
 Note the characterization
 \beql{479bb}
 \frac{1}{2} \ddL=\{ x \in
[0,1]: \dN_j(x) > 0 \textrm{ for all } j \geq 1\},
\eeq
which shows  that $\frac{1}{2} \ddL$  is  a subset of $\ddL$.

The following  weak increasing property will be
used  in the proof of Theorem \ref{th15a} to describe the
points of increase of the Takagi singular function, and in
the proof of Theorem \ref{th37} determining the  expected number of
local level sets per level.

%
%

\begin{theorem}\label{th48a}
Let $x$ belong to the deficient digit set  $ \ddL$.

(1) If $x$ has a binary expansion that does not end in $0^{\infty}$,
then there exists a strictly increasing sequence $\{ x_k\}_{k=1}^{\infty} \subset \ddL$
such that
\begin{equation*}
\lim_{k \to \infty} x_k = x~~~\mbox{and all}~~~ \tau(x_k) < \tau (x).
\end{equation*}

(2) If $x$ has a binary expansion that does not end in $(01)^{\infty}$,
then there exists a strictly decreasing sequence 
$\{ x_k\}_{k=1}^{\infty} \subset \ddL$
such that
\begin{equation*}
\lim_{k \to \infty} x_k = x~~~\mbox{and all}~~~ \tau(x_k) > \tau (x).
\end{equation*}

 (3) If $x \in \frac{1}{2}\ddL$, then (1), (2)  above 
  hold  with the stronger property that all $\{ x_k\}_{k=1}^{\infty} \subset \frac{1}{2} \ddL$.
\end{theorem}

\paragraph{Proof.} 
(1) The condition that the  binary expansion of $x \in \ddL$ not end in $0^{\infty}$
is necessary for there to exist an infinite sequence $x_1 < x_2 < x_3< \cdots \subset \ddL$
such that $\lim_{k \to \infty} x_k = x.$

Write $x=0.b_1 b_2 b_3...$, and note  that there must be infinitely many indices 
$m_1< m_2 < \cdots$ with $\dN_{m_k}(x)>0$.
We then choose the $x_k$ to be the dyadic rationals obtained by suitably
truncating the binary expansion of $x$ at these points:
\begin{equation}
\label{480a}
x_k := 0. b_1 b_2 \cdots b_{m_k}0^{\infty}.
\end{equation}
We clearly have $x_k \in \ddL$ and $\lim_{k \to \infty} x_k =x$, and 
$x_1 \le  x_2 \le x_3 < \cdots$.
This sequence contains an infinite strictly increasing subsequence
since the binary expansion of $x$ does not end in $0^{\infty}$. 
We also note  that if $x \in \frac{1}{2} \ddL,$ i.e. if each $\dN_j(x) \ge 1,$
then each $x_k \in \frac{1}{2} \ddL.$ 

It remains  to show that $\tau(x_k) < \tau(x).$ By  Lemma \ref{le21} we have 
\beql{482}
\tau(x) = \sum_{j=1}^{\infty} \frac{\ell_j}{2^j}.
\eeq
Letting $N^1(j) := \dM_j(x)$ (resp. $N^{0}(j): = N_j^{0}(x)$) count the number of $1$'s (resp. $0$'s)
in the first $j$ digits of the binary expansion of $x$, we have
\beql{483}
\tau(x_k) = \sum_{j=1}^{m_k} \frac{\ell_j}{2^j} + 
\sum_{j= m_k+1}^{\infty} \frac{N^1(m_k)}{2^j}.
\eeq
Now for all $j > m_k$ we have
$$
N^1(m_k) = \min (N^0(m_k) , N^1(m_k)) \le \min( N^0(j-1), N^1(j-1)) \le \ell_j,
$$
and strict inequality holds here for at least one $j > m_k$ because 
the binary expansion of $x$ does not end in $0^{\infty}$.
We conclude that $\tau(x_k) < \tau (x)$ by
comparing \eqn{482} and \eqn{483} term by term.

\smallskip
(2) The condition that the  binary expansion of $x \in \ddL$ not end in $(01)^{\infty}$
is necessary for there to exist an infinite sequence $x_1 > x_2 >  x_3> \cdots \subset \ddL$
such that $\lim_{k \to \infty} x_k = x.$

Writing $x=0.b_1 b_2 b_3 \cdots $, there must be infinitely many indices
$m_1 < m_2 < m_3 < \cdots$ such that
\beql{485}
\dN_{j}(x) \ge \dN_{m_k}(x) ~~~\mbox{for all}~~ j \ge m_k.
\eeq
We now choose $x_k$ to be the rational numbers:
\beql{485a}
x_k := 0. b_1 b_2 \cdots b_{m_k} (01)^{\infty}.
\eeq
We clearly have $x_k \in \ddL$ and $\lim_{k \to \infty} x_k =x$, and 
we have $x_1 \ge x_2 \ge x_3 \ge \cdots$
using the fact that $x \in \ddL$ together with \eqn{485}, which
implies that $\dN_j(x_{k+1}) \ge \dN_j(x_k)$ for all $j \ge 1$.  
This sequence contains an infinite strictly decreasing subsequence
since the binary expansion of $x$ does not end in $(01)^{\infty}$. 
Again note that if $x \in \frac{1}{2} \ddL$, so all $\dN_j(x) \ge 1$,
 then all $x_k \in \frac{1}{2} \ddL.$

It remains to show that $\tau(x_k) > \tau(x)$. 
For any $j \ge 1$, set $x[j]:= 0. b_{j+1} b_{j+2} \cdots$ and note that all $x[m_k] \in \ddL$
by virtue of condition \eqn{485}.  
Similarly define $x_k[j]$ and note that 
$x_k[m_k] =0.(01)^{\infty} \in \ddL$.
There are  two cases.

\smallskip
{\em Case (i).} If $\dN(m_k)=0$ then $m_k = 2m$ and 
$$
\tau(x) = \sum_{j=1}^{2m} \frac{\ell_j}{2^j} + \frac{m}{2^{m_k}} + \frac{1}{2^{m_k}} \tau(x[m_k])
$$
while
$$
\tau(x_k) =  \sum_{j=1}^{2m} \frac{\ell_j}{2^j} + \frac{m}{2^{m_k}} +
\frac{1}{2^{m_k}} \cdot \frac{2}{3}.
$$
Since $\tau(x) \le \frac{2}{3}$ and the only $x \in \ddL$ with $\tau(x)= \frac{2}{3}$
is $x= \frac{1}{3} = 0.(01)^{\infty}$, we conclude that the strict inequality $\tau(x_k) > \tau(x)$
holds. 

\smallskip
{\em Case (ii).} If $\dN(m_k) \ge 1$, then since the first $m_k+1$ digits of $x$ and $x_k$ match
we have, using Lemma \ref{le21}, 
$$
\tau(x) = \sum_{j=1}^{m_k} \frac{\ell_j}{2^j} +  \frac{N^1(m_k)}{2^{m_k}}+ 
\frac{\dN_{m_k}(x)}{2^{m_k}} {x[m_k]}
+  \frac{1}{2^{m_k} } \tau(x[m_k]).
$$
and
$$
\tau(x_k) = \sum_{j=1}^{m_k} \frac{\ell_j}{2^j} +  \frac{N^1(m_k)}{2^{m_k}}+
  \frac{\dN_{m_k}(x)}{2^{m_k}} x_k[m_k] +
 \frac{1}{2^{m_k}}  \tau( x_k[m_k]).
 $$
 Now $x_k[m_k]= 0.(01)^{\infty} = \frac{1}{3} > x[m_k]$ and $\tau(x_k[m_k]) = \frac{2}{3} \ge
 \tau(x[m_k])$, so we conclude the strict inequality $\tau(x_k) > \tau(x)$, as required.
 
 (3) Suppose $x \in \frac{1}{2}\ddL.$ 
 For (1) any truncation $x_k$ given
 by \eqn{480a}
 will automatically satisfy the defining property \eqn{479bb} 
 for membership in  $\frac{1}{2} \ddL.$
 Similarly for (2) any value $x_k$ given by  \eqn{485a} will automatically satisfy \eqn{479bb}.
 $~~~\Box$\\

Next consider the Takagi function restricted to the set
 $\frac{1}{2} \ddL.$
  We show  that  the Takagi function is  nondecreasing on this set,
  and moreover is strictly increasing off a certain specific countable
  set of $x$.
We thank P. Allaart for the following  proof  to establish the
nondecreasing property,  which replaces our original  argument.

%
%

\begin{theorem}\label{th49}
  (1) The Takagi function  is nondecreasing on the set $\frac{1}{2} \ddL$.
  
   (2) The Takagi function  is strictly increasing on $\frac{1}{2} \ddL$
  away from a countable set of points, which are a subset of  those rationals having binary
  expansions ending in $0^{\infty}$ or $(01)^{\infty}$. For each level
  $y$ the  equation $y= \tau(x)$
   has at most two solutions with $x \in \frac{1}{2} \ddL$. Thus  if $x_1< x_2< x_3 $ are
  all in $\frac{1}{2}\ddL$ then $\tau(x_3) > \tau(x_1)$.

\end{theorem}

\paragraph{Proof.}\hspace{-.5em}
(1) [Allaart] \hspace{.5em}
Let $x < x' \in \frac 1 2 \ddL$ have the binary expansions 
$x=0.b_1b_2\ldots$ and $x' = 0.b_1'b_2'\ldots$,  and let $n$ be the index of the first bit that
differs in the two expansions: that is, $b_i = b_i'$ for $i < n$ and
$b_n = 0$, $b_n' = 1$. 
Now set $x_1 := 0.b_1'\ldots b_n'= \frac{k+1}{2^n},$  for some $k \ge 0$. Clearly, 
$x< x_1 \leq x'$.  Furthermore, 
$\dN_i(x_1)= \dN_i(x')> 0$ for $1 \le i \le n$, while
$\dN_i(x_1)= \dN_n(x_1)+(i-n) >0$ for
all $i >n$, so that 
$x_1 \in \frac{1}{2}\ddL$.  Lemma \ref{le23} applied to the interval
$[\frac {k+1}{2^n}, \frac{k+2}{2^n}]$ gives $\tau_n(x_1) \leq
\tau_n(x')$ since $D_n(x_1) = D_n(x') > 0$; hence
$$\tau(x_1) = \tau_n(x_1) \leq \tau_n(x') \leq \tau(x').$$
 Therefore to prove the 
nondecreasing property, it is enough to prove $\tau(x) \le \tau(x_1)$.
We prove the following stronger claim (which does not require 
that either $x$ or $x_1 \in \frac{1}{2} \ddL$).

\paragraph{Claim.}
If  $\frac {k} {2^n} \leq x <\frac{k+1}{2^n}$ and $\dN_i(x) > 0$ for all $i \geq n$, then $\tau(x)
\leq \tau(\frac {k+1}{2^n})$.  

\paragraph{Proof of claim.}
The  result is proved by induction on the value of $m=\dN_n(x)$, at each
step proving it for all $n \ge m.$ We use the self-affine formula
of Lemma \ref{le25}, on $[\frac{k}{2^n}, \frac{k+1}{2^n}]$. Setting $x_0=\frac{k}{2^n}$,
for   $x= x_0+ \frac{w}{2^n}$ with $0 \le w \le 1$,  it gives
$
\tau(x) = \tau(x_0) + \frac{1}{2^n}( \tau(w) + D_n(x_0) w).
$
Taking $w=1$ gives
$\tau(\frac{k+1}{2^n})= \tau(x_0) +\frac{1}{2^n} D_n(x_0),$
and differencing yields
\beql{491}
\tau(\frac{k+1}{2^n}) - \tau(x) = \frac{1}{2^n} \left(D_n(x_0) -(\tau(w) + D_n(x_0) w)\right).
\eeq
Here we note that
\beql{492}
\dN_{i}(w) = \dN_{n+i}(x) -\dN_{n}(x_0), ~~~ ~~ i \ge 1.
\eeq
That is, the function $\dN(\cdot)$ itself undergoes a linear  shift under the
variable change from $x$ to $w$.

We begin the induction with the base case
$m = 1$. We then have $m=\dN_{n}(x) = \dN_{n}(x_0) = 1$,
so that for any $n \ge 1$, \eqn{491} becomes
\beql{493}
\tau(\frac{k+1}{2^n}) - \tau(x)=\frac{1}{2^n} \left(1 -(\tau(w) + w)\right).
\eeq
Using \eqn{492} the   assumption $\dN_i(x) >0$ for $i \geq n$ implies that
$\dN_{i}(w) \ge 0$ for all $i \ge 1$. This says that $w \in \ddL$, so we have $w \le \frac{1}{3}$
whence  $\tau(w) + w \le \frac{2}{3} + \frac{1}{3} = 1. $ Substituting this  inequality in \eqn{493}
completes the base case.

For the inductive step, fix $m \geq 1$ and assume the claim holds
for all $\dN_n(x) = m$ and all $n \ge m$.
Now suppose $\dN_n(x) = m+1$.  We bisect the line segment 
$\left [ \frac {k}{2^n}, \frac {k+1} {2^n}\right)
= \left[ \frac {k}{2^n}, \frac{2k +1}{2^{n+1}} \right) \cup \left[ \frac{2k + 1}{2^{n+1}},
  \frac{k+1}{2^n}\right)$ into the sections where the function
  $f_n(x) := \tau_{n+1}(x) -  \tau_{n}(x) = \frac { \ll 2^n x \gg}{2^n}$ has constant slopes $+1$
  and $-1$ respectively and check the claim in these two cases:  
  
  {\em Case (i).}
   Suppose $ \frac{2k+1}{2^{n+1}} \le x < \frac{k+1}{2^{n}}.$
  Here  $\dN_{n+1}(x) = m$, since $\dN_{j}(x)$ is the slope
  of $\tau_{j}(x)$ at $x$ (see Lemma \ref{le23}).  The claim assumption
  gives  $\dN_{i}(x) > 0$ for all $i \ge n+1$, so the 
  induction hypothesis now applies to $x$ at level $n'=n+1,$ to give
   $\tau(x) \leq \tau(\frac{2k+2}{2^{n+1}})=  \tau(\frac{k+1}{2^n})$.  
  
  {\em Case (ii).}
  Suppose  $\frac{k}{2^n} \le  x <  \frac{2k+1}{2^{n+1}}$. Now \eqn{491}  becomes
  $$
  \tau(\frac{k+1}{2^n}) - \tau(x)=\frac{1}{2^n} \Large(m+1 -(\tau(w) + (m+1)w)\Large),
  $$
  and the Case (ii) range of $x$ implies $0 \le w \le  \frac{1}{2}.$ 
  However for this range of $w$ we have
   $$
  \tau(w) + (m+1)w \le 1+ \frac{m+1}{2} \le m+1.
    $$
  Substituting this in the previous inequality gives $\tau(x) \le \tau(\frac{k+1}{2^n})$.
 (Note: no conditions on $\dN_{i}(x)$ are required in Case (ii).)
  This completes the induction step, and the claim follows.\smallskip

  (2) By Theorem \ref{th48a}(3) we have $\tau(x_1) < \tau(x_2)$
  for any two points $x_1$ and $x_2$ in $\frac{1}{2} \ddL$ with
  $x_1 < x_2$ such that
     neither $x_1$ nor $x_2$ is a rational number with binary expansion
   ending in either 
   $0^{\infty}$ or $(01)^{\infty}$. 
   More is true; the conditions in Theorem \ref{th48a}
    imply furthermore  that equality $t=\tau(x_1) = \tau(x_2)$ for $x_1< x_2$
   in $\frac{1}{2} \ddL$ can occur only if $x_1$ ends in $(01)^{\infty}$
   and $x_2$ ends in $0^{\infty}$. Thirdly, using the
   nondecreasing property of $\tau(x)$ on $\frac{1}{2} \ddL$, 
   we infer that  for any level $y$ the equation
    $y= \tau(x)$ has at most two solutions $x \in \frac{1}{2} \ddL$.
   (A countable set of values $y$ having two solutions  in $\frac{1}{2}\ddL$ exists, 
   with solutions being pairs
  $(\frac{1}{2}x(B)^{-}, \frac{1}{2} x(B)^{+})$ associated to all intervals $I(B), B \in  \sB$.)
   $~~~\Box$

%
%
%
%
\section{ Takagi Singular Function}\label{sec5}
\setcounter{equation}{0}

We now study the behavior of the Takagi function restricted to the left hand
endpoints of all local level sets. This leads to defining a singular function whose
points of increase are confined to these endpoints, 
which we name the {\em Takagi singular function.}


\subsection{Flattened Takagi function and Takagi singular function} \label{sec51}

We consider the Takagi function restricted to the set $\ddL$ and linearly
interpolate it across all intervals removed in constructing $\ddL$, obtaining
a new function, the  flattened Takagi function, as follows.

%
%

\begin{defi}~\label{de34}
{\em
The 
{\em $\ddL$-projection function} $P^{L}(x) := x_b$ where $x_b$
is the largest point $x_b \in \ddL$ having $x_b \le x$. 
This function is well-defined since $\ddL$ is a closed set, and $0 \in \ddL$.
It clearly has the projection property
$$
P^{L}( P^{L}(x)) = P^L(x). 
$$
}
\end{defi}

To compute $P^{L}(x)$,
if  $x= 0.b_0b_1 ...$ has $x \not\in \ddL$  then we have
 $$
  x_b = 0.b_0b_1... b_n(01)^\infty,
 $$
  where $n$ is the smallest location
in the binary expansion of $x$ 
such that 
$\dN_j(x) \ge 0$ for $j \le n$, but $\dN_{n+1}(x) <0$.
If no such $n$ exists then $x_b = x \in \ddL$.

%
%

\begin{defi} \label{de35}
{\em 
The {\em flattened Takagi function} $\tauL(x)$ is given by
$$
\tauL(x) := \tau(x_b) -(x-x_b) = \tau (P^{L}(x))  +P^{L}(x)- x.
$$
}
\end{defi}

\begin{figure}[h]
  \begin{center}                                                               
   $ \begin{array}{cc}
           \includegraphics[width=2.0in]{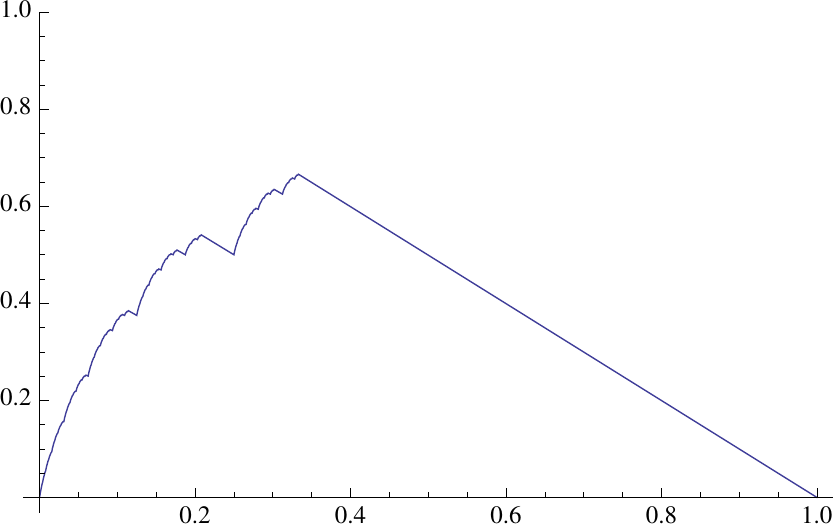} &
           \includegraphics[width=2.0in]{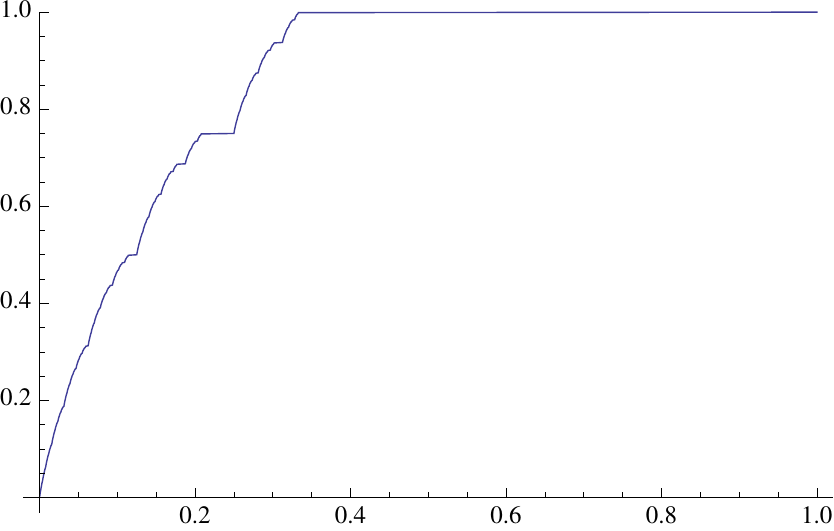} 
         \end{array}$
  \end{center}
  \caption{Graph of flattened Takagi function $\tauL(x)$ (left) and Takagi singular function 
 $\tauS(x)$  (right) .}
  \label{fig411}
\end{figure}

This definition agrees with the Takagi function on the set
$\ddL$ of left hand endpoints of local level sets, 
and it is linear with  slope $-1$ across all omitted intervals $I_{B}$ between
such endpoints. According to
Lemma \ref{le32}(2) this function then
  linearly interpolates across those intervals, 
 showing that $\tauL(x)$ is a continuous function.
It captures all the variation of the Takagi function that is
outside local level sets and 
``flattens it "  inside local level sets. 
It is pictured in Figure \ref{fig411}. (See also the data
in Table \ref{tab2}.)

%
%

 The flattened Takagi function has a countable set of
 intervals on which it has slope $-1$. It seems natural to adjust
 the function to have slope $0$ on these intervals. 
 Thus we study the following function.
 
%
%

\begin{defi} \label{de36}
{\em 
The {\em  Takagi singular function}  $\tauS(x)$ is given by 
\begin{equation*}
\tauS(x) := \tauL(x) + x,
\end{equation*}
where $\tauL$ is the flattened Takagi function.
That is, $\tauS(x) = \tau(x_b) + x_b$.}
\end{defi} 
It is also pictured in Figure \ref{fig411}.
\smallskip

%
%

We now  prove Theorem \ref{th15a} in the introduction,
showing that $\tauS(x)$ is indeed a
singular  continuous function,
justifying its name. 
Note that the  statement of Theorem \ref{th15a} gave a different definition of 
$\tauS$,  and part of the proof below is to show that this alternate definition coincides
with the function given by Definition \ref{de36}.
The distributional derivative  of  $\tauS$ is a  singular measure $\muT$,  the Takagi singular measure, 
which we  study 
in \cite{LM10b}.

\paragraph{Proof of Theorem \ref{th15a}.}
We first show that the function $\tauS(x)$ given by Definition \ref{de36}
is a monotone singular
function.
Lemma~\ref{le32}~(2) implies
that the  flattened Takagi function has slope $-1$ across all omitted intervals outside
of $\ddL$, hence  the definition of the singular Takagi function given in
Definition \ref{de36} guarantees that
it is linear with slope $0$ across all such intervals  $I_B$ for   $B \in \sB$.
These intervals were shown to
have full Lebesgue measure in Theorem~\ref{th33}~(3), so its variation is
confined to the set $\ddL$, which is of measure $0$.  To conclude it is
a Cantor function it  remains to show that
$\tauS(x)$ is a nondecreasing function on $[0,1]$.  Since it is
constant away from $\ddL$, it suffices to show that $\tauS(x)$ is
nondecreasing when restricted to $\ddL$.  Now when $x \in \ddL$
we have $\tauS(x) = \tau(x) + x$, and
by the dyadic self-similarity equation
in Lemma~\ref{le22}, 
$\tau(x)+ x = 2\tau(\frac{x}{2})$, and here $\frac{x}{2} \in \frac{1}{2} \ddL$.
Therefore 
 the nondecreasing property of $\tauS(x)$ is equivalent
 to showing
that $\tau(x)$ is nondecreasing when restricted to $x \in \frac{1}{2}\ddL$.
This was  shown in  Theorem \ref{th49}.

Next we show  that the function $\tauS(x)$ given by Definition \ref{de36}
coincides with the function defined in the theorem statement. That is, it 
has the two properties: (1) $\tauS(x) = \tau(x) + x$ for all $x \in \ddL$;
and (2)
$
\tauS(x) = \sup\{ \tauS(x_1):  x_1 \le x ~~~\mbox{with}  ~~ x_1 \in \ddL\}.
$
The first property holds since for $x \in \ddL$
one has $x=x_b$ so  $\tauS(x) = \tau(x_b) + x_b= \tau(x)+x.$
The second property holds since
$\tauS(x)$ is now known to be nondecreasing, whence
$$
\tauS(x) = \tau(x_b) + x_b 
=\sup\{ \tauS(x_1):  x_1 \le x ~~~\mbox{with}  ~~ x_1 \in \ddL\}.
$$

Finally we verify that the set $\ddL$ is the closure of the set of points
of increase of $\tauS(x)$. It follows from Theorem \ref{th48a} 
that all points of $\ddL$ are points of increase
except for a countable set of 
$x \in \ddL$ which are 
rational numbers  whose binary expansion ends in $0^{\infty}$
or $(01)^{\infty}$.  
Since $\ddL$ is a perfect set, these rational 
numbers are limit points of elements
of $\ddL$ that are irrational, hence they fall in the closure of the set of points of increase
of $\ddL$. $~~~\Box$.

  \paragraph{\bf Remark.} An alternate proof of  the 
  monotone property of $\tauS(x)$ in  Theorem~\ref{th15a} can
  be based on defining piecewise linear approximation functions
  $\tauL_{n}(x)$ and  $\tauS_n(x)$ to $\tauL(x)$ and 
  $\tauS(x)$, respectively, in an obvious fashion. One can prove by induction
  that each $\tauS_{n}(x)$ is a nondecreasing function,  using Lemma \ref{le25}(1).
  The approximations $\tauS_{n}(x)$ approach $\tauS(x)$ pointwise from below,
  giving the result.
  
%
%

\subsection{Jordan decomposition of flattened Takagi function} \label{sec53}

We prove that the flattened Takagi function is of bounded
variation.
Recall that
a  {\em function of bounded (pointwise) variation}
 $f$ on $U= (0,1)$ is a (possibly discontinuous) function whose
 {\em total variation}, denoted ${\rm Var} \, f $ or $V_{0}^{1}(f)$, given by 
\begin{equation*}
V_{0}^1(f) = {\rm Var} \, f:= \sup \{ \sum_{i=1}^n |f(x_i)- f(x_{i-1})| :  ~0<x_0 < x_1 < \cdots < x_n<1, n \ge 1 \}
\end{equation*}
is finite.  
We let $BPV((0,1))$ denote the set of  functions of bounded
(pointwise)  variation on the open interval 
$(0,1)$, following the notation of Leoni \cite[Chap 2]{Le09}.
(In the literature this space is usually denoted $BV(I)$, but this
notation leads to a conflict
with the geometric measure theory notation in Section \ref{sec6N}.)
Any function of bounded variation has a {\em monotone decomposition}
(or {\em Jordan decomposition}) 
\begin{equation*}
f= f_u + f_d
\end{equation*}
in which $f_u$ is an upward monotone (i.e. non-decreasing) bounded function,
possibly with jump discontinuities,  and $f_d$ is a downward monotone
(i.e. non-increasing) bounded function.  Such a decomposition is not unique.
Conversely, any function having a Jordan decomposition is of bounded pointwise variation.
A {\em minimal monotone decomposition} is
one such that 
$$
V_0^1(f) = V_{0}^1(f_u) + V_{0}^1(f_d).
$$

%
%
\begin{theorem}\label{th36} {\em (Jordan decomposition of flattened Takagi function)}
The flattened Takagi function $\tauL(x)$ is of bounded (pointwise) variation,
so is in $BPV((0,1)).$
It has a minimal monotone decomposition  given by
\beql{344}
\tauL(x) = f_u(x) + f_d(x),
\eeq
with downward part
$f_d(x)= -x$ and upward part
 $f_u(x)= \tau^L(x)+x$
both being continuous functions. The upward part is a singular function
whose points of increase are supported on the deficient digit set $\ddL$.
The total
variation of the flattened Takagi function is $V_0^1(\tauL)= 2,$ with $V_0^1(f_u(x))= V_0^1(f_d(x))=1$.
\end{theorem}

%
%
\paragraph{Proof.}
(1) The decomposition \eqn{344} holds by definition of $\tauS(x)$. 
By Theorem \ref{th15a}
  $f_u(x)$ is non-decreasing and 
bounded, and is a monotone singular function supported on $\ddL$.
Clearly $f_d(x) = -x$ is  non-increasing and bounded, 
thus $\tauL(x)$ is of bounded variation, hence \eqn{344} is a monotone 
decomposition. 

(2) The minimality of the monotone decomposition \eqn{344}  is a consequence of the fact
that the function  $f_d$ is 
absolutely continuous with respect to Lebesgue measure,
while by Theorem \ref{th15a}
 the function $f_u$ is singular with respect to  Lebesgue measure.
 Thus
$V_0^1(f) = V_0^1(f_u) + V_0^1(f_d) = 1+1= 2,$
as asserted. $~~~\Box$\\

%
%

\section{Expected number of local level sets} \label{sec6N}
\setcounter{equation}{0}

The results of the last section show that the Takagi function
restricted to the set $\ddL$ is  well behaved,
giving a function of bounded variation.
We  now use the coarea formula 
of geometric measure theory for BV-functions to determine
the expected number of local level sets on a random level 
$0 \le y \le \frac{2}{3}.$

We formulate the coarea formula 
for  the one-dimensional case, using the terminology of 
Leoni \cite{Le09}; note that  the full power of this formula  lies in
the $n$-dimensional case. 
For an open set $U$ of
the real line, the bounded variation space in the sense of geometric
measure theory $BV(U)$ consists of
all functions $f \in L^{1}(U)$ for which there exists a finite signed
Radon measure $\mu$ on the Borel sets of $U$ such that
$$
\int_{U} f(x) \phi'(x) dx = - \int_{U} \phi(x) d\mu
$$
for all test functions $\phi(x) \in C_{c}^{1}(U)$, where we let $C_c^k(U)$
denote the set of $k$-continuously differentiable functions with
compact support.  This signed measure
$\mu$ 
is called the {\em weak derivative} of $f(x)$ and is denoted
$Df$.

The Jordan decomposition theorem for measures
says that the  finite  Radon measure
$Df$  decomposes uniquely into the difference  $Df= Df^{+} - Df^{-}$
of mutually singular nonnegative measures,
 both of which are finite (\cite[Theorem B.72]{Le09}).  The associated 
{\em total variation measure} $|Df|$ is $|Df|= Df^{+} + Df^{-}.$ (See  
Evans and Gariepy \cite[Sect 5.1, Theorem 1]{EG92} for the
$n$-dimensional case.)
The total variation of a function $f$  in $BV(U)$ is expressible 
using test functions as
\begin{equation}\label{test-function-definition}
|Df|(U)= V(f, U) := \sup \Big\{ \int_{U} f(x) \phi^{'}(x) dx:
 \phi(x) \in C_c^{1}(U; \RR),
 ||\phi|| = \max_{x \in U}( |\phi(x)|) \le 1\Big\}.
\end{equation}

We define the  {\em perimeter} of a Lebesgue measurable set $E \subset U$,
denoted $|\partial E|(U)$ or $P(E, U)$, to be  the total variation of its characteristic function $\chi_E$ in $U$, i.e.
$|\partial E|(U) = |D \chi_E|(U).$
It follows  from \eqref{test-function-definition} that if
$E$ is an open set in $U= (0,1)$ consisting of a finite number of
non-adjacent intervals (where non-adjacent means no two intervals have a common endpoint),
then the perimeter $|\partial E|(U)$  counts the number
of endpoints of the intervals  inside $U$. (The perimeter  does not detect the endpoints of $U$, e.g. 
the perimeter of $E= (0, \frac{1}{2})$ has $|\partial E|=1$.)
For  general open sets $E \subset U$, which may have infinitely many open intervals,
the value of the perimeter is more complicated.
Two extreme cases are: (1) The  complement $E :=U \smallsetminus C$ of 
the middle-third Cantor set $C$ has perimeter $|\partial E|(U)=0$;
(2) An open set  $E$ having a partition $E = \bigcup_{i=1}^{\infty} U_i$ in which
each open interval $U_i= (a_i, b_i)$ has an adjacent open interval $J_i = (b_i, b_i +\epsilon_i)$ 
disjoint from $E$ necessarily has perimeter $|\partial E|(U)=+\infty$.

Since the classical bounded variation space $BPV(U)$ (c.f. Section
\ref{sec53}) 
consists of functions
while the geometric measure theory space $BV(U)$ 
consists of
equivalence classes of functions agreeing on sets of 
full Lebesgue measure,
these are distinct spaces. 
However they are closely related, as follows
(Leoni  \cite[Theorem 7.2]{Le09}).

%
%

\begin{prop}\label{pr61a} {\em(Relation of $BPV(U)$ and $BV(U)$.)}\\
Let $U \subset \RR$ be an open set.
 (1) If $f: U \to \RR$ is an integrable function belonging to $BPV(U)$, then its
$L^{1}$-equivalence class belongs to $BV(U)$, and satisfies
\begin{equation*}
 \var f \ge |Df|(U).
\end{equation*}

(2) Conversely
any $f \in BV(U)$ has a right continuous representative function $\bar{f}$
in its $L^{1}$-equivalence class that belongs to  $BPV(U)$, and it
satisfies
\begin{equation*}
\var \bar{f} = |Df|(U).
\end{equation*}
\end{prop}

The following is a one-dimensional version of the coarea formula
for functions in $BV(U)$.

%
%

\begin{prop}\label{pr46}
{ \em (Coarea formula for BV functions)}
Let $U=(0,1)$. If $f \in BV (U)$, there holds:

(1) The upper set
\begin{equation*}
E_t:= E_t(f) = \{ x \in U: ~f(x) >t\}
\end{equation*}
has finite perimeter $|\partial E_t|$ for all but a Lebesgue measure $0$ set
of $t \in \RR$, and the mapping
$$
t \mapsto |\partial E_t|(U),    ~~~t \in \RR
$$
is  a Lebesgue measurable function.

(2) In addition the variation measure $|Df|$ of $f$ satisfies
\begin{equation*}
|Df|(U) = \int_{-\infty}^{\infty} |\partial E_t|(U) dt.
\end{equation*}

(3) Conversely, if $f \in L^{1}(U)$ and
$
 \int_{-\infty}^{\infty} |\partial E_t|(U) dt < \infty,
 $
 then $f \in BV(U)$.
\end{prop}

\paragraph{Proof.} Versions of the coarea formula for functions
in BV(U) for a given  open set $U$ in $\RR^n$  are proved in Evans and Gariepy
\cite[Theorem 1, Sec. 5.5]{EG92} and Leoni \cite[Theorem 13.25]{Le09}.
Here we specialize to the case $U=(0,1)$.
$~~~\Box$ \\

In our application, all 
relevant functions  $f \in BPV(U)$ for $U=(0,1)$ are continuous on $[0,1]$,
in which case by Proposition \ref{pr61a}(2) we have
$|Df|(U) = \var f= V_{0}^1(f)$. 

We  use the Coarea formula for BV functions to
 compute the  expected number of local level sets 
at a random level over the range $0 \le y \le \frac{2}{3}$, with respect to Lebesgue measure.

%
%
\begin{theorem}\label{th37} {\em (Expected number of local level sets)}
For a full Lebesgue measure set of ordinate points $y \in [0, \frac{2}{3}]$ the
number $N^{loc}(y)$ of local level sets at level $y$ is
finite. Furthermore $N^{loc}(y)$ is a Lebesgue measurable function satisfying
\beql{361}
\int_{0}^{\frac{2}{3}}  N^{loc}(y) dy = 1.
\eeq
That is,  the expected number of local level sets  on 
 a randomly drawn ordinate level $y$  is $\frac{3}{2}$.
\end{theorem}

\paragraph{Proof.} 

Theorem \ref{th36} shows
that  the  flattened Takagi function $\tauL(x)$  belongs to $BPV(U)$,
for $U=(0,1).$  Using Proposition \ref{pr61a}(1) we may view it in $BV(U)$ and 
 apply  the  coarea formula for BV functions (Proposition \ref{pr46}), taking 
 $f= \tauL$ and $U= (0,1)$
to obtain
\beql{365}
 \int_{-\infty}^{\infty} |\partial E_t (\tau^L)|(U) dt=|D\tauL|(U),
\eeq
in which $|\partial E_t (\tau^L)|(U)$ is a Lebesgue measurable
function of $t$.
Since the function $\tauL(x)$ is continuous on $[0,1]$,
Proposition \ref{pr61a}(2) applies to give
$$
|D\tauL|(U) = V_{0}^{1}(\tauL)= 2.
$$
Thus we obtain
\begin{equation*}
 \int_{0}^{\frac{2}{3}}  |\partial E_t (\tau^L)|(U) dt =2.
\end{equation*}

We now study to what extent the integrand $|\partial E_t(\tauL)|(U)$ detects 
endpoints of local level sets at level $t$. The  function  
$N^{\loc}(t)$ takes nonnegative integer values or $+\infty$.

\smallskip
{\bf Claim.} {\em For irrational $t$ with  $0< t < \frac{2}{3}$, there holds
\begin{equation*}
2N^{\loc} (t) = |\partial E_t(\tau^L)|(U).
\end{equation*}
 In particular, $N^{loc}(t)$ is a Lebesgue
measurable function, because it differs from 
 $\frac{1}{2} |\partial E_t(\tau^L)|(U)$ on a set of measure zero.}
\smallskip

To prove the claim, we note that for $t >0$ the upper set $E_t := E_t(\tauL)= \bigcup_{i} I_i$ is a finite or
countable disjoint union of nonempty open
intervals, which all lie strictly inside $U=(0,1)$.
Let $I=(a, b)$ be one such interval. 
If $a \notin \ddL$, then the point $(a, \tauL(a))$ 
would lie in the interior of a line segment of slope $-1$ in
the graph $\sG(\tauL)$.  This is impossible since a point in the interior of a line 
of slope $-1$ in $\sG(\tauL)$ will have $\tauL(a+ \epsilon) = t- \epsilon < t$ for 
all small enough positive $\epsilon$, contradicting $I \subseteq
E_t$.  Hence $a \in \ddL$, and consequently $a$ is the left hand
endpoint of a local 
level set in the level set $\tau(x) = t$.
We next assert that  the point $b$ must correspond to a point in the
interior of a line segment of slope $-1$ in the graph $\sG(\tauL)$. 
This holds because, first,  it cannot
correspond to an endpoint of such a segment because then it would be a rational number
with expansion ending in $0^{\infty}$, contradicting  $\tau(b)=t$ being  irrational,
and second,  it cannot be a member of
$\ddL$, since otherwise Theorem~\ref{th48a} (1) would imply there exists
an increasing sequence of values $\{x_k\} \subseteq \ddL$ with $x_k \to b$,
having $\tau(x_k)< \tau(b)$, which would contradict  $b$ being the right
endpoint of an interval of an upper set.
Thus we know  $b$ corresponds to  a point in the interior of a segment of slope $-1$ in $\sG(\tau_L)$, 
hence there is some positive  $\epsilon$ depending on $b$  such that the
adjacent  open interval $(b, b+ \epsilon)$  is not contained in  $E_t$.
We conclude immediately that if
$E_t$ is an infinite disjoint union of nonempty open intervals,
then $|\partial E_t|(U) = \infty$. Moreover $N^{\loc}(t) =
\infty$, proving the claim in this case.

To finish, we may assume $E_t$ is a finite disjoint union of nonempty
open intervals $E_t = \cup_{i = 1}^n (a_i, b_i)$, and hence its boundary $\partial
E_t :=  \bar{E}_t \smallsetminus E_t= \{a_1,\ldots, a_n, b_1,\ldots, b_n\}$.  By the above, each 
$a_i \in \ddL$, but each $b_j \notin \ddL$, so there are gaps
between each of these intervals, and we conclude the perimeter
$|\partial E_t|(U) = 2n$. 
It remains to  show $N^{loc}(t)=n$. The fact that
each $a_i \in \ddL$ implies $N^{loc}(t) \ge n$.
Now let $x \in \ddL$ be the left hand
endpoint of any local level set in the level set $\tau(x) = t$.  
It suffices to  show that $x = a_i$ for some $1 \leq i \leq n$ since then
$N^{\loc}(t) \le \# \{a_1,\ldots, a_n\} = n$. 
We know  $x$ is irrational because $\tau(x) = t$ is irrational by hypothesis.  
Theorem~\ref{th48a}(2) then applies to show there
is a decreasing sequence $\{x_k\} \subseteq \ddL$ with
$x_k \to x$ and $\tau(x_k) > \tau(x) = t$.    
The existence of the sequence $x_k$ implies that $x \in \partial E_t = \{a_1,\ldots, a_n, b_1,\ldots, b_n\}$ since 
$x \notin E_t$ but is given as a limit from of decreasing values $x_k \in E_t$.  By the above, $x \neq b_j$ for any $j$, 
and therefore $x = a_i$ for some $1 \leq i \leq n$, finishing the proof of the claim.

The claim together with \eqn{365} gives \eqn{361}. $~~~\Box$.

\paragraph{Remark.}
Theorem~\ref{th37} gives no information
concerning  the multiplicity of local level sets on those levels having an uncountable level
set, because the set of such levels $y$ has Lebesgue measure $0$.

%
%
%
%

\section{Levels Containing Infinitely Many Local Level Sets}\label{sec5B}
\setcounter{equation}{0}

In this section we show 
  that there exists a dense set of levels in $[0, \frac{2}{3}]$
  that contain a  countably infinite
number of local level sets; this complements Theorem~\ref{th37}.
 We first show this holds for the  particular  level set $L(\frac{1}{2})$.
 The fact that this level set is countably infinite  was previously
 noted by  Knuth \cite[Sec. 7.2.1.3, Problem 82e]{Kn09}.
%
%

\begin{theorem}~\label{th38} {\em (Countably infinite level set)}
The level set  $L(\frac{1}{2})$ is 
countably infinite, with $L(\frac{1}{2}) = \sL_1 \cup \sL_2$ and
\begin{equation*}
\sL_1:= \left\{ x_k :=\frac{1}{2} - \sum_{j=1}^k (\frac{1}{4})^j :  k= 0, 1, 2, ..., \infty \right\},
\end{equation*}
with $\sL_2 = \{ 1-x: ~x \in \sL_1\}.$ 
It contains an infinite number of distinct local
level sets.  One has $L_{x_{\infty}}^{loc}=\{\frac{1}{6} , \frac{5}{6}\} \subset L(\frac{1}{2})$ 
and $\frac{1}{6} \le x \le \frac{5}{6}$ for all $x \in L(\frac{1}{2}).$

\end{theorem}

\paragraph{Proof.} First we show that each $x \in \sL_1 \cup \sL_2$
satisfies $\tau(x)= 1/2$.  By Lemma \ref{le22}(1), it suffices to consider $x \in \sL_1$.  
Let $x := x_k := \frac 1 2 - \sum_{j=1}^k (\frac 1 4)^j$ for $0 \leq k \le  \infty$. 
For $k \ge 1$ the dyadic rational $x_k$ has two binary expansions, the
first expansion  
being  $x_k^{+}=0 .0 (01)^{k-1} 1$ and the second being
$x_k^{-}=0.0(01)^{k-1}0 1^{\infty}$. 
Here the
first binary expansion clearly has all $\dN_j(x_k^{+}) \ge 0$ for $j \ge 1$ 
which  certifies  that $x_k^{+} \in \ddL$. By direct calculation
$\tau(x_0)= \tau(x_1) =\frac{1}{2}$. Now the flip operation shows 
for $k \ge 1$ that 
$$
x_{k}^{-} =0.0(01)^{k-1} 01^{\infty} =0.(0 (01)^k 1 )1^{\infty}  
\sim 0. 0(01)^k 1 0^{\infty} = x_{k+1}^{+}.
$$
Thus we deduce
$$
\tau(x_{k}^{+}) = \tau( x_k^{-}) = \tau(x_{k+1}^{+}),
$$
whence by induction on $k \ge 1$ we conclude that 
$\tau(x_k) = \tau(x_k^{+})= \tau(x_k^{-})=\frac{1}{2}$ for all finite
$k \ge 1$, as asserted. Finally the case $k=\infty$, with $x_{\infty}=0.0(01)^{\infty} = \frac{1}{6}$
has 
$\tau(x_{\infty}) = \lim_{k \to \infty} \tau(x_k)= \frac{1}{2}$ by continuity of the Takagi function. 

Next we observe that the local level sets are
 $L_{x_\infty}^{loc} = \{\frac 1 6, \frac 5 6\}$ and 
$$L_{x_k^{+}}^{\loc} = \{ x_k^+, \, x_{k-1}^-, \, 1- x_k^+, \,
1-x_{k-1}^- \},~~\textrm{ for }k \geq 1.$$
For example, when $k=1$ this is 
$L_{x_1^{+}}^{\loc} = \{\left(\frac{1}{4}\right)^+, 
\left(\frac{1}{2}\right)^-, \left(\frac{3}{4}\right)^-, \left(\frac 1
2 \right)^+ \}.$
This shows that $L(\frac{1}{2})$ has infinitely many
local level sets.

It remains to  show that $L(\frac 1 2)$ contains no elements other than those in
$\sL_1 \cup \sL_2$.   For this, in view of the symmetry of $\tau(x)$, it suffices to prove
two assertions:

\begin{enumerate}
\item[(1)] If $x < \frac 1 6 = x_\infty$,
  then $\tau(x) < \frac 1 2$.
\item[(2)] For all $k$, if $x_k < x < x_{k-1}$, then $\tau(x) > \frac 1 2$.
\end{enumerate}

Let $x < \frac 1 6$.  Lemma \ref{le22}(1) implies that 
$ \tau( x ) = \frac 1 2 \tau(2x) + x$.
Combining this with the inequality $ \tau(x) \leq \frac 2 3$ proves (1).

If $x$ satisfies $x_k < x < x_{k-1} = x_k + \frac {1}{2^{2k}}$ for $k \geq
1$, then $x= x_k + \frac {x'}{2^{2k}}$ for $0 < x' < 1$.    Since
$\dN_{2k}(x_k) = 0$, then by Lemma \ref{le25},
$$\tau(x) = \tau(x_k) + \frac{\tau(x')}{2^{2k}} > \frac{1}{2}.$$
This proves (2).$~~~\Box$\\

Now define
\beql{621}
\Lambda_{\infty}^{loc} := \{ y:  L(y) ~\mbox{contains infinitely many different local level sets} \}.
\eeq
Theorem \ref{th38} above  shows that $y= \frac{1}{2} \in
\Lambda_{\infty}^{loc}.$
Also, recall from Definition \ref{de32} that the breakpoint set $\sBp$
is the set of all balanced dyadic rationals $B' \in \ddL$.
%
%

\begin{theorem}~\label{th62a} {\em (Levels with an infinite number of  local level sets)}

(1) The set $\Lambda_{\infty}^{loc}$ 
has Lebesgue measure $0$. It is not a closed set. 

(2) For each ${B'}=0.b_1 b_2 \cdots b_{2m}\in \sBp$,
 the value
\begin{equation*}
y_{B'} :=  \tau({B'}) + \frac{1}{2^{2m+1}}
\end{equation*}
is a dyadic rational, and $L(y_{B'})$ contains
infinitely many disjoint local level sets. Furthermore, there are
infinitely many   dyadic rationals $x \in \ddL$ with $y_{B'}=\tau(x).$

(3) The set of levels 
\begin{equation*}
\Delta_{\infty}^{\loc} := \left\{ y_{B'}= \tau({B'}) + \frac{1}{2^{2m+1}}:  {B'} \in \sBp  \right\}
\end{equation*}
is dense in $[0, \frac{2}{3}]$. Since $\Delta_{\infty}^{loc}  \subseteq  \Lambda_{\infty}^{loc} $,
the set $\Lambda_{\infty}^{loc}$ is dense in $[0, \frac{2}{3}]$.
\end{theorem}

\paragraph{Proof.}
(1) This measure $0$ property of $\Lambda_{\infty}^{loc}$ follows  immediately from the expected number of local
level sets being finite (Theorem~\ref{th37}). 
The fact that this set is not a closed set will follow once property (3) is proved.\smallskip

(2) For each balanced dyadic rational ${B'}=0.b_1 b_2 \cdots b_{2m}$ in  $\ddL$ we 
consider for $k \ge 1$ the infinite set of dyadic rationals
$$
x_k({B'}):  =0.b_1 b_2 \cdots b_{2m} 0(01)^{k-1} 1 = {B'} + \frac{1}{2^{2m}} x_k,
$$
where $x_k= 0.0(01)^{k-1}1$ has $\tau(x_k)= \frac{1}{2}$ by Theorem~\ref{th38}. 
Using the self-affine scaling property in Lemma \ref{le25}  we have
$$
\tau(x_k({B'})) = \tau({B'}) + \frac{1}{2^{2m}} \tau(x_k) = \tau({B'}) + \frac{1}{2^{2m+1}},
$$
so all points $\tau(x_k({B'}))$ are on the same level $y= y_{B'},$ and $y_{B'}$ is necessarily a
 dyadic rational number. 
Clearly each $x_k({B'}) \in \ddL$, so each determines a different local level set,
establishing (2).\smallskip

(3) It is easy to see that the  set of balanced dyadic rationals ${B'}=0. b_1... b_{2m}$ 
having $\dN_j({B'}) \geq 0$ for all $j \ge 0$ and $\dN_{2m}=0$ 
is dense inside the deficient digit set $\ddL$. Indeed, given 
any $x=0.b_1 b_2... \in \ddL$, the approximation $x_k= 0.b_1 b_2 ... b_k 1^{\dN_k(x)} 0^{\infty}$
is such a dyadic rational having $|x-x_k| \le 2^{-k}.$
Since there
is at least one local level set on each level, we have $\tau(\ddL) = [0, \frac{2}{3}].$
Since the flattened Takagi function is continuous, we conclude that 
the values $y_{B'} = \tau({B'} +\frac{1}{2^{2m+1}})$ are dense in $[0, \frac{2}{3}],$ as asserted.
$~~~\Box$\\

\paragraph{Remarks.}

\smallskip
(1) The proof above shows the stronger result that  if $y \in \Lambda_{\infty}^{loc}$ then
for every balanced dyadic rational ${B'}$ that belongs to $\ddL$, one has 
$y_{B'}^{\ast} : = {B'} + \frac{y}{2^{2m}}  \in \Lambda_{\infty}^{loc}.$ 

\smallskip
 (2) One can 
ask whether the equality $\Delta_{\infty}^{loc} = \Lambda_{\infty}^{loc}$ might hold,
or (weaker) whether $\Lambda_{\infty}^{loc}$ is a countable set.

%
%
%
%

\section{Further Questions}\label{sec9}
\setcounter{equation}{0}

This investigation of the  structure of local level sets of the Takagi function raises a number of   questions for further work.

\medskip
(1) Theorem \ref{th15} shows that an abscissa generic local level set is
uncountable with probability one, with $x$ drawn uniformly from $[0,1]$.
Can one determine the expected number of local level sets at a 
level $L(\tau(x))$, with $x$ drawn uniformly from $[0,1]$?

\medskip
We note that Theorem \ref{th16} sheds no light regarding this question.
It seems to involve properties  related to  a new 
measure $\nu$,  supported on $\ddL$, which is  mutually singular to both
Lebesgue measure and to
the Takagi singular measure $\muT$, which we hope to discuss
elsewhere.

\medskip
(2) Theorem~\ref{th16} shows that the expected number of local level sets
at a given height $y$ drawn uniformly in $[0, \frac{2}{3}]$  is $\frac{3}{2}$.
There is an associated probability distribution 
$$
\Prob[ N^{loc}(y) = k] := \frac{3}{2} \meas[ y: N^{loc}(y) =k ],
$$
whose mean value is $\frac{3}{2}$. Can one explicitly compute these
probabilities in closed form?

\medskip
(3) Can one explicitly determine the Hausdorff dimension of the local level set $L_x^{loc}$
in terms of properties of the binary expansion of $x$? In particular,  to 
what extent  does the balance set $Z(x)$ determine the Hausdorff dimension of
$L_x^{loc}$? 

\medskip
For rational $x$, we have calculated the Hausdorff dimension of
$L_x^{loc}$ in
the proof of Theorem~\ref{th50}.  For general $x$,
recall  that a  necessary condition for positive Hausdorff dimension given in
the proof of Theorem~\ref{th15} is
that $\limsup_{k \to \infty} \frac{k}{c_k} > 0$; this condition
depends only on $Z(x)$.
\medskip

 \medskip
 (4) Theorem \ref{th50} characterizes those rationals $x$ which have
 an uncountable local 
 level set $L^{\loc}$ in terms of their binary expansions. 
 Can one explicitly characterize (e.g. in terms of binary expansion)
 which rational levels $y$ contain some
 rational $x$ for which $L_x^{\loc}$ is 
 uncountable?

\medskip
This problem, which is
a weaker version of one proposed by Knuth \cite[Sect. 7.2.1.3, Exercise 83]{Kn09},  
may be difficult. 

\medskip
(5) The Fourier series of the Takagi function, viewed as a periodic function of
period $1$,  is explicitly known in closed form.
Can one explicitly find the Fourier series of the flattened Takagi function $\tauL(s)$
or the Takagi singular function $\tauS(x)$?

\medskip
The structure and behavior 
of monotone singular functions, particularly including
their Fourier transforms,  is a topic of some interest, tracing back to
 work of  Hartman and Kershner \cite{HK37}, 
 Salem \cite{Sal42}, \cite{Sal43}, \cite{Sal51}. See Dovgoshey et al  
 \cite{DMRV06} for a detailed treatment of the Cantor function.

\paragraph{Acknowledgments.}
We thank  D. E. Knuth for raising  questions  on
the  Takagi function to one of us. We thank S. T. Kuroda for remarks on Takagi's
original construction, and
Mario Bonk for helpful remarks on Hausdorff dimension
and BV functions. We thank  Pieter Allaart for  allowing us to include
his simplified proof of Theorem \ref{th49},  and for bringing the work 
of Buczolich to our attention. We thank the two reviewers
for many insightful comments and corrections, and for some
additional references.

%
%
%
%


\noindent 
$$
\begin{array}{lll}
\mbox{Jeffrey C. Lagarias}    & ~& \mbox{Zachary Maddock}\\
\mbox{Dept. of Mathematics} &~& \mbox{Dept. of Mathematics} \\
\mbox{The University of Michigan} &~& \mbox{Columbia University}\\
\mbox{Ann Arbor, MI 48109-1043} &~& \mbox{New York, NY 10027}\\
\mbox{email: {\tt lagarias@umich.edu}} &~& \mbox{email: {\tt maddockz@math.columbia.edu}}
\end{array}
$$

\end{document}